# Predictor-Based Output Feedback for Nonlinear Delay Systems

**Iasson Karafyllis[*] and Miroslav Krstic[**]**

[*]Dept. of Environmental Eng., Technical University of Crete,
73100, Chania, Greece, email: ikarafyl@enveng.tuc.gr

[**]Dept. of Mechanical and Aerospace Eng., University of California,
San Diego, La Jolla, CA 92093-0411, U.S.A., email: krstic@ucsd.edu

**Abstract**

We provide two solutions to the heretofore open problem of stabilization of systems with arbitrarily long delays at the input and output of a nonlinear system using output feedback only. Both of our solutions are global, employ the predictor approach over the period that combines the input and output delays, address nonlinear systems with sampled measurements and with control applied using a zero-order hold, and require that the sampling/holding periods be sufficiently short, though not necessarily constant. Our first approach considers general nonlinear systems for which the solution map is available explicitly and whose one-sample-period predictor-based discrete-time model allows state reconstruction, in a finite number of steps, from the past values of inputs and output measurements. Our second approach considers a class of globally Lipschitz strict-feedback systems with disturbances and employs an appropriately constructed successive approximation of the predictor map, a high-gain sampled-data observer, and a linear stabilizing feedback for the delay-free system. We specialize the second approach to linear systems, where the predictor is available explicitly. We provide two illustrative examples—one analytical for the first approach and one numerical for the second approach.

## 1. Introduction

<u>Summary of Results of the Paper.</u> Even though numerous results have been developed in recent years for stabilization of nonlinear systems with input delays by *state feedback* [18,20,24,25,26,27,30,31,32,46,50], and although additional delays in state measurements are allowed in our recent work [20], the problem of stabilization of systems with arbitrarily long delays at the input and/or output by *output feedback* has remained open.

We provide two solutions to this problem. Both of our solutions address nonlinear systems with sampled measurements and with control applied using a zero-order hold, with a requirement that the sampling/holding periods be sufficiently short, though not necessarily constant. Both of our solutions also employ the predictor approach to provide the control law with an estimate of the future state over a period that combines the input and output delays.

Our first approach considers general nonlinear systems for which the solution map is available explicitly and whose one-sample-period predictor-based discrete-time model allows state reconstruction, in a finite number of steps, from the past values of inputs and output measurements.

Our second approach considers a class of globally Lipschitz strict-feedback systems with disturbances and employs an appropriately constructed successive approximation of the predictor map, a high-gain sampled-data observer, and a linear stabilizing feedback for the delay-free system. The results of the second approach can be applied to the linear time-invariant case as well, providing robust global exponential sampled-data stabilizers, which are completely insensitive to perturbations of the sampling schedule.

Both of our approaches achieve global asymptotic stabilization. The first approach also achieves dead-beat stabilization in case the delay-free plant is dead-beat stabilizable. The second approach achieves input-to-state stability with respect to plant disturbances and measurement disturbances, as well as global exponential stability in the absence of disturbances.

<u>Problem Statement and Literature.</u> As in [18,20,24,25,26,27,50] we consider nonlinear systems of the form:

$$\dot{x}(t) = f(x(t), u(t-\tau))$$
$$x(t) \in \Re^n, u(t) \in U \subseteq \Re^m$$
(1.1)

where $\tau \geq 0$ is the input delay, $U \subseteq \Re^m$ is a non-empty set with $0 \in U$ (the control set) and $f : \Re^n \times U \to \Re^n$ is a locally Lipschitz mapping with $f(0,0) = 0$. We employ the predictor-based approach, which is ubiquitous for linear systems (see [40] and the references in [25,26]) and is different from other approaches for systems with input delays [30,31,32,46], where the stabilizing feedback for the delay free system is either applied or is modified and stability is guaranteed for sufficiently small input delays. The input in (1.1) can be applied continuously or with zero-order hold (see [20]) and the measured output is usually assumed to be the state vector $x(t) \in \Re^n$.

In [20], we extended predictor-based nonlinear control to the case of sampled measurements and measurement delays expressed as

$$y(t) = x(\tau_i - r), \text{ for } t \in [\tau_i, \tau_{i+1})$$

where $y$ is the measured output, the discrete time instants $\tau_i$ are the sampling times and $r \geq 0$ is the measurement delay. The motivation is that sampling arises simultaneously with input and output delays in control over networks. Few papers have studied this problem (exceptions are [14] where input and measurement delays are considered for linear systems but the measurement is not sampled and [22] where the unicycle is studied).

In the absence of delays, in sampled-data control of nonlinear systems semiglobal practical stability is generally guaranteed [10,35,36,37], with the desired region of attraction achieved by sufficiently fast sampling. Alternatively, global results are achieved under restrictive conditions on the structure of the system [9,13,39]. Simultaneous consideration to sampling and delays (either physical or sampling-induced) is given in the literature on control of linear and nonlinear systems over networks [7,8,12,37,39,44,45,49], but almost all available results rely on delay-dependent conditions for the existence of stabilizing feedback. Exceptions are the papers [3,28], where prediction-based control methodologies are employed.

The assumption that the state vector is measured is seldom realistic. Instead, measurement is a function of the state vector, i.e., the measured output of system (1.1) is given by:

$$y(t) = h(x(iT_1 - r)), \quad t \in [iT_1, (i+1)T_1), i \in Z^+$$
(1.2)

where $T_1 > 0$ is the sampling period, $r \geq 0$ is the measurement delay and $h: \Re^n \to \Re^k$ is a continuous vector field with $h(0) = 0$ (the output map). Notice that the measurements are obtained at discrete time instants.

We study the following problem in this paper: find a feedback law, which utilizes the sampled measurements and applies the input with zero-order hold, given by

$$u(t) = u_j, \ t \in [jT_2, (j+1)T_2), j \in Z^+ \tag{1.3}$$

where $T_2 > 0$ is the holding period, such that the closed-loop system (1.1) with (1.2), (1.3) is globally asymptotically stable.

Two Solutions Provided in the Paper. The above problem is solved for two particular cases:

1$^{st}$ Case (Section 2): The case where the solution map of the open-loop system (1.1) is *explicitly known* and $T_1 = T_2 = T > 0$ (sample-and-hold case). Under appropriate assumptions for observability, forward completeness and sampled-data stabilizability of the open loop system (1.1) with $\tau = 0$, we guarantee stabilization of system (1.1) with a predictor-based version of any sampled-data controller designed for the delay-free plant. For example, all sampled-data feedback designs proposed in [9,10,13,21,35,36,37,39] which guarantee global stabilization can be exploited for the stabilization of a delayed system with input/measurement delays, sampled measurements and input applied with zero order hold. The class of feedforward systems (see [23,26] and references therein) can be addressed by using the proposed observer-based predictor feedback design.

2$^{nd}$ Case (Sections 3,4 and 5): The class of globally Lipschitz systems of the form

$$\begin{aligned} \dot{x}_i(t) &= f_i(x_1(t),...,x_i(t)) + x_{i+1}(t) + g_i(x(t),u(t))d_i(t) \ , \ i = 1,...,n-1 \\ \dot{x}_n(t) &= f_n(x(t)) + g_n(x(t),u(t))d_n(t) + u(t-\tau) \\ x(t) &= (x_1(t),...,x_n(t))' \in \Re^n, u(t) \in \Re, d(t) = (d_1(t),...,d_n(t))' \in \Re^n \end{aligned} \tag{1.4}$$

where $f_i: \Re^i \to \Re$ ($i=1,...,n$) are globally Lipschitz functions with $f_i(0) = 0$ ($i=1,...,n$) and the output map is $h(x) = x_1$. The inputs $d_i$ ($i=1,...,n$) represent disturbances and the functions $g_i: \Re^i \to \Re$ ($i=1,...,n$) are locally Lipschitz, bounded functions. In this case, we can show stabilizability of system (1.1) even under arbitrary perturbations of the sampling schedule, by combining the sampled-data observer design in [17] and the approximate predictor control proposed in [18]. We also show robustness with respect to measurement errors and modeling errors. The feedback design is based on the corresponding delay free system

$$\begin{aligned} \dot{x}_i(t) &= f_i(x_1(t),...,x_i(t)) + x_{i+1}(t) \ , \ i = 1,...,n-1 \\ \dot{x}_n(t) &= f_n(x(t)) + u(t) \end{aligned} \tag{1.5}$$

The proposed control schemes for both cases consist of three components:

1$^{st}$ Component: An observer, which utilizes past input and output values in order to provide (continuous or discrete) estimates of the delayed state vector $x(t-r)$.

2$^{nd}$ Component: The predictor mapping that utilizes the estimation provided by the observer and past input values in order to provide an estimation of the future value of the state vector $x(t+\tau)$.

*3rd Component:* A nominal globally stabilizing feedback for the corresponding delay-free system.

The above control scheme has long been in use for linear systems [29,33,34,48,51] and it has been used even for partial differential equation systems [11], but is novel for nonlinear systems. Moreover, even for Linear Time-Invariant (LTI) systems

$$\dot{x}(t) = Ax(t) + Bu(t-\tau) + Gd(t) \qquad (1.6)$$

where $x(t) \in \Re^n, u(t) \in \Re, d(t) \in \Re^n$, we provide new sampled-data feedback exponential stabilizers that are robust to perturbations of the sampling schedule.

*Notation.* Throughout the paper we adopt the following notation:

* For a vector $x \in \Re^n$ we denote by $|x|$ its usual Euclidean norm, by $x'$ its transpose. For a real matrix $A \in \Re^{n \times m}$, $A' \in \Re^{m \times n}$ denotes its transpose and $|A| := \sup\{|Ax|; x \in \Re^n, |x|=1\}$ is its induced norm. $I \in \Re^{n \times n}$ denotes the identity matrix. By $A = \text{diag}(l_1, l_2, ..., l_n)$ we mean a diagonal matrix with $l_1, l_2, ..., l_n$ on its diagonal.

* $\Re_+$ denotes the set of non-negative real numbers. $Z^+$ denotes the set of non-negative integers. For every $t \geq 0$, $[t]$ denotes the integer part of $t \geq 0$, i.e., the largest integer being less or equal to $t \geq 0$. A partition $\pi = \{T_i\}_{i=0}^{\infty}$ of $\Re_+$ is an increasing sequence of times with $T_0 = 0$ and $T_i \to +\infty$.

* We say that an increasing continuous function $\gamma : \Re_+ \to \Re_+$ is of class $K$ if $\gamma(0) = 0$. We say that a function $\gamma$ of class $K$ is of class $K_\infty$ if $\lim_{s \to +\infty} \gamma(s) = +\infty$. By $KL$ we denote the set of all continuous functions $\sigma : \Re_+ \times \Re_+ \to \Re_+$ with the properties: (i) for each $t \geq 0$ the mapping $\sigma(\cdot, t)$ is of class $K$; (ii) for each $s \geq 0$, the mapping $\sigma(s, \cdot)$ is non-increasing with $\lim_{t \to +\infty} \sigma(s,t) = 0$.

* By $C^j(A)$ ($C^j(A;\Omega)$), where $A \subseteq \Re^n$ ($\Omega \subseteq \Re^m$), $j \geq 0$ is a non-negative integer, we denote the class of functions (taking values in $\Omega \subseteq \Re^m$) that have continuous derivatives of order $j$ on $A \subseteq \Re^n$.

* Let $x:[a-r,b) \to \Re^n$ with $b > a \geq 0$ and $r \geq 0$. By $T_r(t)x$ we denote the "history" of $x$ from $t-r$ to $t$, i.e., $(T_r(t)x)(\theta) := x(t+\theta); \theta \in [-r,0]$, for $t \in [a,b)$. By $\tilde{T}_r(t)x$ we denote the "open history" of $x$ from $t-r$ to $t$, i.e., $(\tilde{T}_r(t)x)(\theta) := x(t+\theta); \theta \in [-r,0)$, for $t \in [a,b)$.

* Let $I \subseteq \Re$ be an interval. By $L^\infty(I;U)$ ($L^\infty_{loc}(I;U)$) we denote the space of measurable and (locally) bounded functions $u(\cdot)$ defined on $I$ and taking values in $U \subseteq \Re^m$. Notice that we do not identify functions in $L^\infty(I;U)$ which differ on a measure zero set. For $x \in L^\infty([-r,0];\Re^n)$ or $x \in L^\infty([-r,0);\Re^n)$ we define $\|x\|_r := \sup_{\theta \in [-r,0]} |x(\theta)|$ or $\|x\|_r := \sup_{\theta \in [-r,0)} |x(\theta)|$. Notice that $\sup_{\theta \in [-r,0]} |x(\theta)|$ is not the essential supremum but the actual supremum and that is why the quantities $\sup_{\theta \in [-r,0]} |x(\theta)|$ and $\sup_{\theta \in [-r,0)} |x(\theta)|$ do not coincide in general.

* The saturation function $sat(x)$ is defined by $sat(x) = x/|x|$ for all $x \in \Re$.

Throughout the paper, for $r = 0$ we adopt the convention $L^\infty([-r,0];\Re^n) = \Re^n$ and $C^0([-r,0];\Re^n) = \Re^n$. Finally, for reader's convenience, we mention the following fact, which is a direct consequence of Lemma 2.2 in [1] and Lemma 3.2 in [15]. The fact is used extensively throughout the paper.

**FACT:** *Suppose that the system $\dot{x}(t) = f(x(t), u(t))$ is forward complete. Then for every $x_0 \in \Re^n$, $u \in L^\infty_{loc}([-\tau,+\infty);\Re^m)$, the solution $x(t)$ of (1.1) with initial condition $x(0) = x_0 \in \Re^n$ and corresponding to input $u \in L^\infty_{loc}([-\tau,+\infty);\Re^m)$ exists for all $t \geq 0$. Moreover, for every $T > 0$ there exists a function $a \in K_\infty$ such that for every $x_0 \in \Re^n$, $u \in L^\infty_{loc}([-\tau,+\infty);\Re^m)$ the solution $x(t)$ of (1.1) with initial condition $x(0) = x_0 \in \Re^n$ and corresponding to input $u \in L^\infty_{loc}([-\tau,+\infty);\Re^m)$ satisfies $|x(t)| \leq a\left(|x_0| + \sup_{-\tau \leq s \leq t-\tau} |u(s)|\right)$, for all $t \in (0,T]$.*

## 2. Solution Map Known Explicitly

We consider system (1.1) under the following hypotheses:

**Hypothesis (H1):** *The system*

$$\dot{x}(t) = f(x(t), u(t)), \quad x(t) \in \Re^n, u(t) \in U \subseteq \Re^m \tag{2.1}$$

*is forward complete.*

**Hypothesis (H2):** *There exists $k : \Re^n \to U$, $a \in K_\infty$ such that the $0 \in \Re^n$ is Globally Asymptotically Stable for the closed-loop system (2.1) with*

$$u(t) = k(x(iT)), \quad t \in [iT,(i+1)T), \quad i \in Z^+ \tag{2.2}$$

*i.e., there exists $\sigma \in KL$ such that for every $t \geq 0$, $x_0 \in \Re^n$, the solution $x(t)$ of the closed-loop system (2.1) with (2.2) at time $t \geq 0$ with initial condition $x(0) = x_0 \in \Re^n$ exists and satisfies $|x(t)| \leq \sigma(|x_0|, t)$. Moreover, the following inequality holds:*

$$|k(x)| \leq a(|x|), \quad \forall x \in \Re^n \tag{2.3}$$

Let $\phi(t, x_0; u)$ denote the solution map of (2.1), i.e., the unique solution $x(t) \in \Re^n$ of (2.1) at time $t \geq 0$ with initial condition $x(0) = x_0 \in \Re^n$ and corresponding to a measurable and essentially bounded input $u : [0,t] \to U$ satisfies $x(t) = \phi(t, x_0; u)$. The control approach that we will use for the stabilization of system (1.1) assumes explicit knowledge of the solution map $\phi(t, x_0; u)$ of (2.1). If the output map were the identity function then the approach in [20] could be directly applied for the stabilization of (1.1). Here, we need an additional observability hypothesis.

Let $l \in Z^+$ be an integer such that $r + \tau = lT + \delta$, where $\delta \in [0,T)$. If $\delta > 0$, then we can define the operator $P : U^2 \to L^\infty([0,T); U)$ by means of the formula

$$(P(u_1, u_2))(t) := u_1, \text{ for } t \in [0, \delta) \text{ and } (P(u_1, u_2))(t) := u_2, \text{ for } t \in [\delta, T)$$

and the mapping $F: \Re^n \times U^2 \to \Re^n$ by means of the equation:

$$F(x, u_1, u_2) := \phi(T, x; P(u_1, u_2)) \qquad (2.4)$$

Notice that the previous definitions in conjunction with the semigroup property for the solution map, imply for all $i \in Z^+$ with $i \geq l+1$:

$$x((i+1)T - r) = F(x(iT - r), u_1, u_2) \qquad (2.5)$$

where $x(t)$ denotes any solution of (1.1) with piecewise constant input that satisfies $u(iT - r - \tau + \theta) = (P(u_1, u_2))(\theta)$ for all $\theta \in [0, T)$.

If $\delta = 0$ then we similarly define the operator $P: U \to L^\infty([0,T);U)$ by means of the formula

$$(P(u))(t) := u, \text{ for } t \in [0, T)$$

and the mapping $F: \Re^n \times U \to \Re^n$ by means of the equation:

$$F(x, u) := \phi(T, x; P(u)) \qquad (2.6)$$

Notice again that the previous definitions in conjunction with the semigroup property for the solution map, give for all $i \in Z^+$ with $i \geq l$:

$$x((i+1)T - r) = F(x(iT - r), u) \qquad (2.7)$$

where $x(t)$ denotes any solution of (1.1) with constant input that satisfies $u((i-l)T + \theta) \equiv u$.

Therefore for every $\delta \in [0, T)$, we can construct an autonomous discrete-time system of the form

$$\begin{aligned} x(i+1) &= F(x(i), v(i)), \quad x(i) \in \Re^n, v(i) \in V \\ y(i) &= h(x(i)) \end{aligned} \qquad (2.8)$$

which is associated with system (1.1), (1.2) and represents a one-sampling-period "predictor system". Notice that $V = U^2$ for the case $\delta > 0$ and $V = U$ for the case $\delta = 0$. The following observability hypothesis is employed in the present work (see also [16]).

**Hypothesis (H3):** *The discrete-time system (2.8) is completely observable, i.e., there exists $p \in Z^+$, $p \geq 1$ and a continuous function $\Psi: V^p \times \Re^{kp} \times \Re^k \to \Re^n$ with $\Psi(0,0) = 0$ such that the solution of the discrete-time system (2.8) with arbitrary initial condition corresponding to arbitrary input satisfies $x(i + p) = \Psi(v(i), \ldots, v(i + p - 1), y(i), \ldots, y(i + p))$ for all $i \in Z^+$.*

**Example 2.1:** We consider the 3-dimensional feedforward system

$$\begin{aligned} \dot{x}_1(t) &= u(t - \tau), \; \dot{x}_2(t) = x_1(t) + x_1(t)u(t - \tau), \; \dot{x}_3(t) = x_2(t) + x_1^2(t) \\ x &= (x_1, x_2, x_3)' \in \Re^3, u \in U \subseteq \Re \end{aligned} \qquad (2.9)$$

This system (for $\tau = 0$) is not feedback linearizable [23]. Hypothesis (H1) holds for system (2.9) and its solution map with $\tau = 0$ is given by:

$$\phi_1(t,x;u) = x_1 + \int_0^t u(s)ds$$

$$\phi_2(t,x;u) = x_2 + t\,x_1 + \int_0^t \int_0^\tau (1+u(\tau))u(s)ds\,d\tau + x_1 \int_0^t u(s)ds$$

$$\phi_3(t,x;u) = x_3 + x_2 t + \frac{t^2}{2}x_1 + \int_0^t \int_0^w \int_0^\tau (1+u(\tau))u(s)ds\,d\tau\,dw + 3x_1 \int_0^t \int_0^\tau u(s)ds\,d\tau + x_1^2 t + \int_0^t \left(\int_0^w u(s)ds\right)^2 dw$$

(2.10)

Here we study the case $r+\tau \in (0,T)$, in which the equality $r+\tau = lT + \delta$ holds with $\delta = r+\tau$, $l = 0$ and the mapping $F:\Re^3 \times \Re^2 \to \Re^3$ defined by (2.4) is given by

$$F_1(x,u_1,u_2) = x_1 + Q_1(u_1,u_2)$$
$$F_2(x,u_1,u_2) = x_2 + T\,x_1 + x_1 Q_1(u_1,u_2) + G_2(u_1,u_2)$$
$$F_3(x,u_1,u_2) = x_3 + T(x_2 + x_1^2) + \frac{T^2}{2}x_1 + 3x_1 Q_2(u_1,u_2) + G_3(u_1,u_2)$$

(2.11)

for all $(x,u_1,u_2) \in \Re^n \times U^2$, where $Q_1, Q_2, G_2, G_3 : \Re^2 \to \Re$ are defined by:

$$Q_1(u_1,u_2) := \delta u_1 + (T-\delta)u_2$$

$$Q_2(u_1,u_2) := \frac{\delta^2}{2}u_1 + \delta(T-\delta)u_1 + \frac{(T-\delta)^2}{2}u_2$$

$$G_2(u_1,u_2) = Q_2(u_1,u_2) + \frac{\delta^2}{2}u_1^2 + \delta(T-\delta)u_1 u_2 + \frac{(T-\delta)^2}{2}u_2^2$$

$$G_3(u_1,u_2) = \frac{\delta}{2}\left[T^2 - T\delta + \frac{\delta^2}{3}\right]u_1 + u_2\frac{(T-\delta)^3}{6} + \frac{3\delta^2}{2}\left(T - \frac{2\delta}{3}\right)u_1^2 + 3u_1 u_2 \delta\frac{(T-\delta)^2}{2} + u_2^2\frac{(T-\delta)^3}{2}$$

We consider two cases:

$1^{st}$ Case (two states are measured): $U \subseteq \Re$ with $0 \in U$ is arbitrary and the measured output of (2.9) is given by (1.2), where

$$y = (y_1, y_2)' = h(x) = (x_1, x_3)'$$

(2.12)

In this case, system (2.8) is completely observable with $p = 1$, since the solution of the discrete-time system (2.8) with arbitrary initial condition and corresponding to arbitrary input satisfies

$$x(i+1) = \begin{bmatrix} y_1(i+1) \\ T^{-1}(y_2(i+1) - y_2(i)) - y_1^2(i) + y_1(i)B(u_1(i),u_2(i)) + C(u_1(i),u_2(i)) \\ y_2(i+1) \end{bmatrix}$$

(2.13)

where

$$B(u_1,u_2) = -3T^{-1}Q_2(u_1,u_2) + Q_1(u_1,u_2) + \frac{T}{2} \quad , \quad C(u_1,u_2) = G_2(u_1,u_2) - T^{-1}G_3(u_1,u_2)$$

$2^{nd}$ Case (only one state is measured): $U \subseteq [-\varepsilon, \varepsilon]$ with $0 \in U$, $\varepsilon \in \left(0, \frac{1}{6}\right)$ and the measured output of (2.9) is given by (1.2), where

$$y = h(x) = x_3$$

(2.14)

In this case the following equation holds for all $i \in Z^+$ for the solution of the discrete-time system (2.8), where $F: \Re^3 \times \Re^2 \to \Re^3$ is defined by (2.11):

$$x_3(i+2) - 2x_3(i+1) + x_3(i) - P(u_1(i), u_2(i), u_1(i+1), u_2(i+1)) =$$
$$\left(T^2 + 3TQ_1(u_1(i), u_2(i)) + 3Q_2(u_1(i+1), u_2(i+1)) - 3Q_2(u_1(i), u_2(i))\right) x_1(i)$$

where $P(u_1, u_2, v_1, v_2) := G_3(v_1, v_2) - G_3(u_1, u_2) + TG_2(u_1, u_2) - \frac{1}{2}\left(T^2 + 6Q_2(v_1, v_2)\right) Q_1(u_1, u_2) - TQ_1^2(u_1, u_2)$.

The above definitions $Q_1(u_1, u_2) := \delta u_1 + (T - \delta) u_2$, $Q_2(u_1, u_2) := \frac{\delta^2}{2} u_1 + \delta(T - \delta) u_1 + \frac{(T-\delta)^2}{2} u_2$ show that if $U \subseteq [-\varepsilon, \varepsilon]$ with $\varepsilon \in \left(0, \frac{1}{6}\right)$ then the inequalities $|Q_2(u_1, u_2)| \leq \frac{T^2}{2} \varepsilon$, $|Q_1(u_1, u_2)| \leq T\varepsilon$ hold for all $u_1, u_2 \in U$. Therefore the inequality $T^2 + 3TQ_1(u_1(i), u_2(i)) + 3Q_2(u_1(i+1), u_2(i+1)) - 3Q_2(u_1(i), u_2(i)) \geq (1 - 6\varepsilon)T^2$ holds for all $u_1(i), u_2(i) \in U$, $u_1(i+1), u_2(i+1) \in U$. In this case, system (2.8) is completely observable with $p = 2$, since the solution of the discrete-time system (2.8) with arbitrary initial condition and corresponding to arbitrary input satisfies

$$x(i+2) = \begin{bmatrix} M + Q_1(u_1(i+1), u_2(i+1)) \\ T^{-1}(y(i+2) - y(i+1)) - M^2 + M B(u_1(i+1), u_2(i+1)) + C(u_1(i+1), u_2(i+1)) \\ y(i+2) \end{bmatrix} \quad (2.15)$$

$$M = \frac{y(i+2) - 2y(i+1) + y(i) - P(u_1(i), u_2(i), u_1(i+1), u_2(i+1))}{D(u_1(i), u_2(i), u_1(i+1), u_2(i+1))} + Q_1(u_1(i), u_2(i))$$

where $D(u_1, u_2, v_1, v_2) := T^2 + 3TQ_1(u_1, u_2) + 3Q_2(v_1, v_2) - 3Q_2(u_1, u_2)$. ◁

Next, we consider again the general system (1.1), (1.2), (1.3) with $T_1 = T_2 = T > 0$. By virtue of hypothesis (H1) the solution of (1.1), (1.2), (1.3) exists for all $t \geq 0$ for arbitrary initial condition $T_r(0)x = x_0 \in C^0([-r, 0]; \Re^n)$, $\breve{T}_\tau(0)u \in L^\infty([-\tau, 0); U)$ and arbitrary sequence of inputs $\{u_i\}_{i=0}^\infty$. The measurements made up to time $t = iT$ are given by $y(t) = y_j = h(x(jT - r))$, $t \in [jT, (j+1)T)$ for $j = 0, \ldots, i$. Notice that, if we denote $x_j = x(jT - r)$ for all $j \in Z^+$ with $jT \geq r + \tau$, then we obtain $x_{j+1} = F(x_j, u_{j-l-1}, u_{j-l})$ for the case $\delta > 0$ and $x_{j+1} = F(x_j, u_{j-l})$ for the case $\delta = 0$, where $F$ is defined by (2.4) or (2.6) and $l \in Z^+$ is the integer such that $r + \tau = lT + \delta$, where $\delta \in [0, T)$. Hypothesis (H3) implies the existence of a continuous function $R: U^{p+1} \times \Re^{kp} \times \Re^k \to \Re^n$ (called the reconstruction mapping) for the case $\delta > 0$ and $R: U^p \times \Re^{kp} \times \Re^k \to \Re^n$ for the case $\delta = 0$ with $R(0,0) = 0$ such that for every $i \in Z^+$ with $iT \geq (p + l + 1)T$, the following equality holds:

$$x_i = x(iT - r) = R(u_{i-p-l-1}, \ldots, u_{i-l-1}, y_{i-p}, \ldots, y_i), \text{ for the case } \delta > 0 \quad (2.16)$$

$$x_i = x(iT - r) = R(u_{i-p-l}, \ldots, u_{i-l-1}, y_{i-p}, \ldots, y_i), \text{ for the case } \delta = 0 \quad (2.17)$$

We are also in a position to define the predictor mapping that correlates $x(iT - r)$ with $x(iT + \tau)$, which is given by

$$\Phi(x, u_{i-l-1}, \ldots, u_{i-1}) := \phi(r + \tau, x; Q(u_{i-l-1}, \ldots, u_{i-1})), \text{ for the case } \delta > 0 \quad (2.18)$$

$$\Phi(x, u_{i-l}, \ldots, u_{i-1}) := \phi(r + \tau, x; Q(u_{i-l}, \ldots, u_{i-1})), \text{ for the case } \delta = 0 \quad (2.19)$$

where

$$(Q(u_{i-l-1},...,u_{i-1}))(t) := u_{i-l-1}, \text{ for } t \in [0,\delta),$$

$$(Q(u_{i-l-1},...,u_{i-1}))(t) := u_j, \ t \in [(j+l-i)T+\delta,(j+l+1-i)T+\delta), \ j=i-l,...,i-1$$

for the case $\delta > 0$ and for the case $\delta = 0$

$$(Q(u_{i-l},...,u_{i-1}))(t) := u_j, \text{ for } t \in [(j+l-i)T,(j+l+1-i)T), \ j=i-l,...,i-1$$

By virtue of (2.16), (2.17), (2.18) and (2.19) the following equalities hold for every $i \in Z^+$ with $i \geq p+l+1$:

$$x(iT+\tau) = \Phi(R(u_{i-p-l-1},...,u_{i-l-1},y_{i-p},...,y_i),u_{i-l-1},...,u_{i-1}), \text{ for the case } \delta > 0 \quad (2.20)$$

$$x(iT+\tau) = \Phi(R(u_{i-p-l},...,u_{i-l-1},y_{i-p},...,y_i),u_{i-l},...,u_{i-1}), \text{ for the case } \delta = 0 \quad (2.21)$$

The computation of the predictor mapping and the reconstruction mapping is straightforward when the solution map of (2.1) is known. The following example illustrates how easily the prediction and reconstruction mappings can be computed.

**Example 2.2:** We return to Example 2.1 and consider the 3-dimensional feedforward system (2.9) with output defined by (2.14), $U \subseteq [-\varepsilon,\varepsilon]$ with $0 \in U$, $\varepsilon \in \left(0,\frac{1}{6}\right)$, for the case $r+\tau \in (0,T)$. In this case the equality $r+\tau = lT+\delta$ holds with $\delta = r+\tau$, $l=0$ and the mapping $F:\Re^3 \times \Re^2 \to \Re^3$ defined by (2.4) is given by (2.11). In order to define the reconstruction mapping $R:U^3 \times \Re^2 \times \Re \to \Re^3$ we simply need to replace $y(i+2), y(i+1), y(i)$ by $y_i, y_{i-1}, y_{i-2}$, respectively, to use equation (2.15) with $u_2(i) = u_1(i+1)$ and to replace $u_2(i+1), u_2(i), u_1(i)$ by $u_{i-1}, u_{i-2}, u_{i-3}$, respectively. Therefore, we obtain:

$$R(u_{i-3},u_{i-2},u_{i-1},y_{i-2},y_{i-1},y_i) := \begin{bmatrix} M + Q_1(u_{i-2},u_{i-1}) \\ T^{-1}(y_i - y_{i-1}) - M^2 + M\,B(u_{i-2},u_{i-1}) + C(u_{i-2},u_{i-1}) \\ y_i \end{bmatrix} \quad (2.22)$$

$$M = \frac{y_i - 2y_{i-1} + y_{i-2} - P(u_{i-3},u_{i-2},u_{i-2},u_{i-1})}{D(u_{i-3},u_{i-2},u_{i-2},u_{i-1})} + Q_1(u_{i-3},u_{i-2})$$

where the mappings $Q_1, B, C, P, D$ have been defined in Example 2.1. The predictor mapping is simply obtained from (2.10):

$$\Phi_1(x,u_{i-1}) = x_1 + \delta u_{i-1}$$

$$\Phi_2(x,u_{i-1}) = x_2 + \delta x_1 + \frac{\delta^2}{2}(1+u_{i-1})u_{i-1} + \delta x_1 u_{i-1} \quad (2.23)$$

$$\Phi_3(x,u_{i-1}) = x_3 + \delta\left(x_2 + x_1^2\right) + \frac{\delta^2}{2}x_1 + \frac{\delta^3}{6}u_{i-1} + \frac{3\delta^2}{2}x_1 u_{i-1} + \frac{\delta^3}{2}u_{i-1}^2$$

In general, the rule to obtain the reconstruction mapping from the mapping $\Psi:V^p \times \Re^{kp} \times \Re^k \to \Re^n$ for the case $\delta > 0$ involved in Hypothesis (H3) is to replace $v(i+j) = (u_1(i+j),u_2(i+j)) \in V$ by $(u_{i-p-l-1+j},u_{i-p-l+j})$ for $j=0,...,p-1$ (notice that $V = U^2$) and $y(i+j)$ by $y_{i-p+j}$ for $j=0,...,p-1$. ◁

In summary, the proposed control scheme consists of three components:

$1^{\text{st}}$ **Component:** A sampled-data observer based on a state-reconstruction mapping $R:U^{p+1} \times \Re^{kp} \times \Re^k \to \Re^n$, given by (2.16) for the discrete-time one-sample-period "predictor system" (2.8). The reconstruction mapping utilizes past input and output values in order to provide an estimate for the delayed state vector $x(iT-r)$.

2$^{nd}$ Component: The predictor mapping $\Phi$, given by (2.20), which utilizes the estimation provided by the reconstruction map and past input values in order to provide an estimation of the future value of the state vector $x(iT+\tau)$.

3$^{rd}$ Component: The nominal globally stabilizing feedback $k:\Re^n \to U$ involved in Hypothesis (H2), which employs the predictor.

We are now ready to state our main result. Its proof is provided in the Appendix.

**Theorem 2.3:** *Let $T>0$, $r,\tau \geq 0$ with $r+\tau>0$ be given and let $l \in Z^+$ and $\delta \in (0,T)$ be such that $r+\tau = lT+\delta$. Moreover, suppose that Hypotheses (H1)-(H3) hold for system (1.1), (1.2) with $T_1=T_2=T$. Then the closed-loop system (1.1), (1.2), (1.3) with*

$$u(t) = k(\Phi(X, u_{i-l-1},...,u_{i-1})), \text{ for } t \in [iT,(i+1)T) \quad (2.24)$$

*where $u_j = u(jT)$ and*

$$X = R(u_{i-p-l-1},...,u_{i-l-1}, y_{i-p},..., y_i) \quad (2.25)$$

*where $y_j = h(x(jT-r))$, is Globally Uniformly Asymptotically Stable, in the sense that there exists a function $\tilde{\sigma} \in KL$ such that for every $(x_0,u_0) \in C^0([-r-pT,0];\Re^n) \times L^\infty([-(p+l+1)T,0);U)$, the solution $(x(t),u(t)) \in \Re^n \times \Re^m$ of the closed-loop system (1.1), (1.2), (1.3), (2.24), (2.25) with $T_1=T_2=T$, initial condition $\breve{T}_{r+\tau}(0)u = u_0 \in L^\infty([-(p+l+1)T,0);U)$, $T_r(0)x = x_0 \in C^0([-r-pT,0];\Re^n)$ satisfies the following inequality for all $t \geq 0$:*

$$\|T_{r+pT}(t)x\|_{r+pT} + \|\breve{T}_{(p+l+1)T}(t)u\|_{(p+l+1)T} \leq \tilde{\sigma}(\|x_0\|_{r+pT} + \|u_0\|_{(p+l+1)T}, t) \quad (2.26)$$

*Finally, if the closed-loop system (2.1), (2.2) satisfies the dead-beat property of order $jT$, where $j \in Z^+$ is positive, i.e., for all $x_0 \in \Re^n$ the solution $x(t)$ of (2.1), (2.2) with initial condition $x(0) = x_0 \in \Re^n$ satisfies $x(t) = 0$ for all $t \geq jT$, then there exists $q \in Z^+$ such that the closed-loop system (1.1), (1.2), (1.3), (2.24), (2.25) satisfies the dead-beat property of order $qT$, i.e., for every $(x_0,u_0) \in C^0([-r-pT,0];\Re^n) \times L^\infty([-(p+l+1)T,0);U)$, the solution $(x(t),u(t)) \in \Re^n \times \Re^m$ of system (1.1), (1.2), (1.3), (2.24), (2.25) with initial condition $\breve{T}_{r+\tau}(0)u = u_0 \in L^\infty([-(p+l+1)T,0);U)$, $T_r(0)x = x_0 \in C^0([-r-pT,0];\Re^n)$ satisfies $x(t) = 0$ for all $t \geq qT$.*

**Remark 2.4:** A very similar statement holds for the case $\delta = 0$. The only thing that needs to be changed is the feedback law (2.24), (2.25), which is replaced by $u(t) = k(\Phi(X, u_{i-l},...,u_{i-1}))$, for $t \in [iT,(i+1)T)$ and $X = R(u_{i-p-l},...,u_{i-l-1}, y_{i-p},..., y_i)$.

**Example 2.5:** We return to Example 2.2 and consider the 3-dimensional feedforward system (2.9) with output defined by (2.14), $U = [-\varepsilon, \varepsilon]$ with $\varepsilon \in \left(0, \frac{1}{6}\right)$, for the case $r+\tau \in (0,T)$. In this case the equality $r+\tau = lT+\delta$ holds with $\delta = r+\tau$, $l=0$. It is shown in [21] (Theorem 3.7 and Remark 3.8) that Hypothesis (H2) holds for system (2.9) for every $\varepsilon > 0$. More specifically, for every $\varepsilon > 0$ there exist constants $K_0, K_1, K_2, R_1, R_2, T > 0$ such that hypothesis (H2) holds with

$$k(x) := \begin{cases} -K_0 \text{sat}(x_1) & \text{if } |x_1| \geq R_1 \text{ and } \sqrt{x_2^2 + (x_1 + x_2)^2} \geq R_2 \\ -x_1 - K_1 \text{sat}(x_2 + x_1) & \text{if } |x_1| < R_1 \text{ and } \sqrt{x_2^2 + (x_1 + x_2)^2} \geq R_2 \\ -2(x_1 + x_2) - K_2 \text{sat}\left(x_3 + x_2 + \frac{1}{2}x_1\right) & \text{if } \sqrt{x_2^2 + (x_1 + x_2)^2} < R_2 \end{cases} \quad (2.27)$$

Moreover, the inequality $|k(x)| \leq \varepsilon$ holds for all $x \in \Re^3$. It follows from Theorem 2.3 that the closed-loop system (2.9), (1.2), (2.24), (2.25), where $h, R, \Phi$ are defined by (2.14), (2.22), (2.23), respectively, is Globally Uniformly Asymptotically Stable in the sense described in Theorem 2.3.

For the case where $[-\varepsilon, \varepsilon] \subseteq U$ and the output map is given by (2.12) (the case where two states are measured), we showed in Example 2.1 that Hypothesis (H3) holds. It follows from Theorem 2.3 that the closed-loop system (2.9), (1.2) with $u(t) = k(\hat{X}), t \in [iT, (i+1)T)$, where $k : \Re^3 \to \Re$ is defined by (2.27), $\hat{X} = (\hat{X}_1, \hat{X}_2, \hat{X}_3)'$,

$$\hat{X}_1 = X_1 + (r + \tau)u_{i-1}$$
$$\hat{X}_2 = X_2 + (r + \tau)X_1 + \frac{(r+\tau)^2}{2}(1 + u_{i-1})u_{i-1} + (r+\tau)u_{i-1}X_1 \quad (2.28)$$
$$\hat{X}_3 = X_3 + (r+\tau)(X_2 + X_1^2) + \frac{(r+\tau)^2}{2}X_1 + \frac{(r+\tau)^3}{6}(1+u_{i-1})u_{i-1} + \frac{3(r+\tau)^2}{2}u_{i-1}X_1 + \frac{(r+\tau)^3}{3}u_{i-1}^2$$

where $u_{i-1} = u((i-1)T)$ and $X = (X_1, X_2, X_3)'$ is defined by

$$X = \begin{bmatrix} y_1(i) \\ T^{-1}(y_2(i) - y_2(i-1)) - y_1^2(i-1) + y_1(i-1)B(u_{i-2}, u_{i-1}) + C(u_{i-2}, u_{i-1}) \\ y_2(i) \end{bmatrix} \quad (2.29)$$

where $y_1(j) = x_1(jT - r)$, $y_2(j) = x_3(jT - r)$, $j = i-1, i$, $u_{i-2} = u((i-2)T)$ and the functions $B, C : \Re^2 \to \Re$ are defined in Example 2.1, is Globally Uniformly Asymptotically Stable. ◁

## 3. Globally Lipschitz Systems

We consider system (1.4) with output

$$y(\tau_i) = x_1(\tau_i - r) + \xi(\tau_i), \; i \in Z^+ \quad (3.1)$$

where $\{\tau_i\}_{i=0}^{\infty}$ is a partition of $\Re^+$ with $\sup_{i \geq 0}(\tau_{i+1} - \tau_i) \leq T_1$. We assume that $r + \tau > 0$, where $r \geq 0$ is the measurement delay and $\tau \geq 0$ is the input delay. The locally bounded input $\xi : \Re_+ \to \Re$ represents the measurement error and the measurable and locally essentially bounded inputs $d_i : \Re_+ \to \Re$ ($i = 1, \ldots, n$) represent disturbances. We assume that there exist constants $L \geq 0$ and $G \geq 0$ such that

$$|f_i(x_1, \ldots, x_i) - f_i(z_1, \ldots, z_i)| \leq L|(x_1 - z_1, \ldots, x_i - z_i)|, \; \forall (x_1, \ldots, x_i) \in \Re^i, \; \forall (z_1, \ldots, z_i) \in \Re^i \quad (3.2)$$

$$|g_i(x,u)| \leq G, \quad \forall (x,u) \in \Re^n \times \Re \tag{3.3}$$

for all $i=1,...,n$. Moreover, $f_i(0)=0$ for all $i=1,...,n$. Define $f(x):=(f_1(x_1),...,f_n(x))' \in \Re^n$, $A=\{a_{i,j}: i,j=1,...,n\} \in \Re^{n \times n}$ with $a_{i,i+1}=1$ for all $i=1,...,n-1$ and $a_{i,j}=0$ if $j \neq i+1$, $b=(0,...,0,1)' \in \Re^n$, $c:=(1,0,...,0)' \in \Re^n$. We notice that inequalities (3.2), (3.3) guarantee that system (1.4) is forward complete, i.e., for every $(x_0, u, d) \in \Re^n \times L^{\infty}_{loc}([-\tau,+\infty);\Re) \times L^{\infty}_{loc}(\Re_+;\Re^n)$ the solution $x(t) \in \Re^n$ of system (1.4) with initial condition $x(0)=x_0 \in \Re^n$ and corresponding to inputs $(u,d) \in L^{\infty}_{loc}([-\tau,+\infty);\Re) \times L^{\infty}_{loc}(\Re_+;\Re^n)$ exists for all $t \geq 0$. Indeed, notice that the function $P(t) = \frac{1}{2}|x(t)|^2$ satisfies $\dot{P}(t) \leq ((n+1)L+3)P(t) + \frac{1}{2}G^2|d(t)|^2 + \frac{1}{2}u^2(t-\tau)$, for almost all $t \geq 0$ for which the solution $x(t) \in \Re^n$ of system (1.4) exists. Integrating the previous differential inequality and using a standard contradiction argument, we conclude that the solution $x(t) \in \Re^n$ of system (1.4) exists for all $t \geq 0$ and satisfies the following estimate for all $t > 0$:

$$|x(t)| \leq \left( |x_0| + \frac{G \sup_{0 \leq s < t}|d(s)| + \sup_{-\tau \leq s < t-\tau}|u(s)|}{\sqrt{(n+1)L+3}} \right) \exp\left( \frac{(n+1)L+3}{2} t \right) \tag{3.4}$$

The proposed observer/predictor-based feedback law consists of three components:

1) A high-gain sampled-data observer for system (1.4), (3.1) which provides an estimate $z(t) \in \Re^n$ of the delayed state vector $x(t-r)$.

2) An approximate predictor, i.e., a mapping that utilizes the applied input values and the estimate $z(t) \in \Re^n$ provided by the observer in order to provide an estimate for $x(t+\tau)$.

3) A stabilizing feedback law for the delay-free system, i.e., system (1.5).

In what follows, we are going to describe the construction of each one of the components. We also assume that the input and measurement delay values $\tau, r \geq 0$ are perfectly known (see Remark 3.7(a) below for the case where the measurement delay is not precisely known).

**1st Component (High-Gain Sample-Data Observer):** Let $p=(p_1,...,p_n)' \in \Re^n$ be a vector such that the matrix $(A+pc') \in \Re^{n \times n}$ is Hurwitz. The existence of a vector $p=(p_1,...,p_n)' \in \Re^n$ is guaranteed, since the pair of matrices $(A,c)$ is observable. The proposed high-gain sample-data observer is of the form:

$$\begin{aligned}
\dot{z}_i(t) &= f_i(z_1(t),...,z_i(t)) + z_{i+1}(t) + \theta^i p_i(c'z(t) - w(t)), \quad i=1,...,n-1 \\
\dot{z}_n(t) &= f_n(z_1(t),...,z_n(t)) + \theta^n p_n(c'z(t) - w(t)) + u(t-r-\tau) \\
\dot{w}(t) &= f_1(z_1(t)) + z_2(t), \quad t \in [\tau_i, \tau_{i+1}) \\
w(\tau_{i+1}) &= y(\tau_{i+1}) = x_1(\tau_{i+1} - r) + \xi(\tau_{i+1}) \\
\tau_{i+1} &= \tau_i + T_1 \exp(-b(\tau_i)), \tau_0 = 0
\end{aligned} \tag{3.5}$$

where $(z(t), w(t)) \in \Re^n \times \Re$, $\theta \geq 1$ is a constant to be chosen sufficiently large by the user and $b: \Re_+ \to \Re_+$ is an arbitrary non-negative locally bounded input that is unknown to the user. Notice that the sampling sequence $\{\tau_i\}_{i=0}^{\infty}$ is an arbitrary partition of $\Re_+$ with $\sup_{i \geq 0}(\tau_{i+1} - \tau_i) \leq T_1$, i.e., the sampling schedule is arbitrary. In order to justify the use of the high-gain sample-data observer (3.5), we notice that system (3.5) is the feedback interconnection of the usual high-gain observer of system (1.4) which estimates $x(t-r)$ instead of $x(t)$ and uses $w(t)$ instead of (the non-available signal) $x_1(t-r)$:

$$\dot{z}_i(t) = f_i(z_1(t),...,z_i(t)) + z_{i+1}(t) + \theta^i g_i(c'z(t) - w(t)), \ i = 1,...,n-1$$
$$\dot{z}_n(t) = f_n(z_1(t),...,z_n(t)) + \theta^n g_n(c'z(t) - w(t)) + u(t-r-\tau)$$

and the inter-sample predictor of (the non-available signal) $x_1(t-r)$:

$$\dot{w}(t) = f_1(z_1(t)) + z_2(t) \ , \ t \in [\tau_i, \tau_{i+1})$$
$$w(\tau_{i+1}) = x_1(\tau_{i+1} - r) + \xi(\tau_{i+1})$$
$$\tau_{i+1} = \tau_i + T_1 \exp(-b(\tau_i)), \tau_0 = 0$$

which utilizes the measurements and predicts the value of $x_1(t-r)$ between two consecutive measurements. Sampled-data observers of this type (which are robust to sampling schedule perturbations) were first proposed in [17] (see also [41,42,43]).

**2$^{nd}$ Component (Approximate Predictor):** Let $u \in L^\infty([0,T];\mathfrak{R})$ be arbitrary and define the operator $P_{T,u} : C^0([0,T];\mathfrak{R}^n) \to C^0([0,T];\mathfrak{R}^n)$ by

$$(P_{T,u}x)(t) = x(0) + \int_0^t (f(x(\tau)) + Ax(\tau) + bu(\tau))d\tau, \text{ for } t \in [0,T]. \quad (3.6)$$

We denote $P_{T,u}^l = \underbrace{P_{T,u} ... P_{T,u}}_{l \text{ times}}$ for every integer $l \geq 1$. We next define the operators $G_T : \mathfrak{R}^n \to C^0([0,T];\mathfrak{R}^n)$, $C_T : C^0([0,T];\mathfrak{R}^n) \to \mathfrak{R}^n$ and $Q_{T,u}^l : \mathfrak{R}^n \to \mathfrak{R}^n$ for $l \geq 1$ by

$$(G_T x_0)(t) = x_0, \text{ for } t \in [0,T] \text{ and } C_T x = x(T) \quad (3.7)$$

$$Q_{T,u}^l = C_T P_{T,u}^l G_T \quad (3.8)$$

We next define the mapping $P_{l,m}^u : \mathfrak{R}^n \to \mathfrak{R}^n$ for arbitrary $u \in L^\infty([0,r+\tau];\mathfrak{R})$. Let $l,m \geq 1$ be integers and $T = \frac{r+\tau}{m}$. We define for all $x \in \mathfrak{R}^n$:

$$P_{l,m}^u x = Q_{T,u_m}^l ... Q_{T,u_1}^l x \quad (3.9)$$

where $u_i(s) = u(s+(i-1)T)$, $i = 1,...,m$ for $s \in [0,T)$. Notice that $u_i \in L^\infty([0,T];\mathfrak{R})$ for $i = 1,...,m$.

The operator $P_{l,m}^u$ is a nonlinear operator which provides an estimate of the value of the state vector of system (1.5) after $r+\tau$ time units when the input $u \in L^\infty([0,r+\tau];\mathfrak{R})$ is applied. The operator is constructed based on the following procedure:

- first, we divide the time interval $[0,r+\tau]$ into $m \geq 1$ subintervals of equal length $T = \frac{r+\tau}{m}$,

- secondly, we apply the method of successive approximations to each one of the subintervals; more specifically we apply $l \geq 1$ successive approximations in order to get an estimate of the value of the state vector at the end of each one of the subintervals.

**Proposition 3.1 (see [18]):** *Let $l,m$ be positive integers with $(nL+1)T<1$, where $T=\dfrac{r+\tau}{m}$. Then there exists a constant $K:=K(m)\geq 0$, independent of $l$, such that for every $u\in L^\infty([0,r+\tau);\Re)$ and $x\in\Re^n$ the following inequality holds:*

$$\left|P_{l,m}^u x-\phi(r+\tau,x;u)\right|\leq K\frac{((nL+1)T)^{l+1}}{1-(nL+1)T}\left(|x|+\sup_{0\leq\tau<r+\tau}|u(\tau)|\right) \quad (3.10)$$

*where $\phi(t,x;u)$ denotes the unique solution of (1.5) at time $t\in[0,r+\tau]$, with initial condition $x\in\Re^n$ and corresponding to input $u\in L^\infty([0,r+\tau);\Re)$.*

Let $\delta_{r+\tau}:L^\infty_{loc}([-r-\tau,+\infty);\Re)\to L^\infty_{loc}([0,+\infty);\Re)$ denote the shift operator defined by

$$(\delta_{r+\tau}u)(t):=u(t-r-\tau), \text{ for } t\geq 0 \quad (3.11)$$

We are now able to define the approximate predictor mapping $\Phi_{l,m}:\Re^n\times L^\infty([-r-\tau,0);\Re)\to\Re^n$ defined by:

$$\Phi_{l,m}(x,u):=P_{l,m}^{\delta_{r+\tau}u}x \quad (3.12)$$

Using (3.2), (3.3), (3.10) and the Gronwall-Bellman lemma, we conclude that the following inequality holds for the solution of (1.4) for all $t\geq r$:

$$\begin{aligned}\left|\Phi_{l,m}(z,\breve{T}_{r+\tau}(t)u)-x(t+\tau)\right|&\leq K\frac{((nL+1)T)^{l+1}}{1-(nL+1)T}\left(|z|+\sup_{t-r-\tau\leq s<t}|u(s)|\right)\\&+\exp((nL+1)(r+\tau))\left((r+\tau)G\sup_{t-r\leq s\leq t+\tau}|d(s)|+|z-x(t-r)|\right)\end{aligned} \quad (3.13)$$

It should be noticed that by virtue of (3.4) and (3.13), we obtain the following inequality for all $(u,z)\in L^\infty([-r-\tau,0);\Re)\times\Re^n$:

$$\left|\Phi_{l,m}(z,u)\right|\leq\Gamma\left(|z|+\sup_{-r-\tau\leq s<0}|u(s)|\right) \quad (3.14)$$

where $\Gamma:=K\dfrac{((nL+1)T)^{l+1}}{1-(nL+1)T}+\exp\left(\dfrac{(n+1)L+3}{2}(r+\tau)\right)$.

**3rd Component (Delay-Free Stabilizing Feedback):** Due to the triangular structure of system (1.4), the results in [47] in conjunction with (3.2), (3.3), imply that there exists $k\in\Re^n$, a symmetric positive definite matrix $P\in\Re^{n\times n}$ and constants $\mu,\gamma>0$ such that

$$x'P(A+bk')x+x'Pf(x)+x'P\text{diag}(g_1(x,u),...,g_n(x,u))d\leq-2\mu x'Px+\gamma|d|^2, \text{ for all }(x,d,u)\in\Re^n\times\Re^n\times\Re \quad (3.15)$$

We are now in a position to construct a stabilizing observer-based predictor feedback. Let $T_2>0$ be the "holding period". The proposed feedback law is given by (3.5) with

$$u(t)=k'\Phi_{l,m}(z(iT_2),\breve{T}_{r+\tau}(iT_2)u), \text{ for } t\in[iT_2,(i+1)T_2) \quad (3.16)$$

**Theorem 3.2:** *Let $Q\in\Re^{n\times n}$ be a symmetric positive definite matrix that satisfies $Q(A+pc')+(A'+cp')Q+2qI\leq 0$ for certain constant $q>0$ and certain $p\in\Re^n$. Let $P\in\Re^{n\times n}$ be a symmetric positive definite matrix that satisfies (3.15) for certain constant $\mu,\gamma>0$ and certain* $k \in \Re^n$. Let $l, m$ be positive integers and $\theta \geq 1$, $T_2 > 0$ and $T_1 > 0$ be constants such that the following inequalities hold:

$$4|Qp|(L+\theta)T_1\sqrt{\frac{|Q|}{a}} < q \tag{3.17}$$

$$\theta \geq \max\left(1, \frac{2|Q|L\sqrt{n}}{q}\right) \tag{3.18}$$

$$\left((nL+1+|k|)\sqrt{\frac{b'Pb}{2K_1}}+\mu\right)|k|\left(T_2 + K\frac{((nL+1)T)^{l+1}}{1-(nL+1)T}\right) < \mu \tag{3.19}$$

where $a > 0$ is a constant satisfying $a|x|^2 \leq x'Qx$ for all $x \in \Re^n$, $0 < K_1 \leq K_2$ are constants satisfying $K_1|x|^2 \leq x'Px \leq K_2|x|^2$ for all $x \in \Re^n$, $K > 0$ is the constant involved in (3.13) and $T = \frac{r+\tau}{m}$.

Then there exist constants $\Theta_i > 0$ ($i=1,...,6$) and $\sigma > 0$ such that for every $(x_0, u_0, z_0, w_0) \in C^0([-r,0];\Re^n) \times L^\infty([-r-\tau,0);\Re) \times \Re^n \times \Re$, $(\xi, b, d) \in L^\infty_{loc}(\Re_+; \Re \times \Re_+ \times \Re^n)$ the solution $(T_r(t)x, \breve{T}_{r+\tau}(t)u, z(t), w(t)) \in C^0([-r,0];\Re^n) \times L^\infty([-r-\tau,0);\Re) \times \Re^n \times \Re$ of the closed-loop system (1.4), (3.5) and (3.16) with initial condition $\breve{T}_{r+\tau}(0)u = u_0 \in L^\infty([-r-\tau,0);\Re)$, $T_r(0)x = x_0 \in C^0([-r,0];\Re^n)$, $(z(0), w(0)) = (z_0, w_0) \in \Re^n \times \Re$ and corresponding to inputs $(\xi, b, d) \in L^\infty_{loc}(\Re_+; \Re \times \Re_+ \times \Re^n)$ satisfies the following inequality for all $t \geq 0$:

$$|z(t)| + |w(t)| + \|T_r(t)x\|_r + \|\breve{T}_{r+\tau}(t)u\|_{r+\tau} \leq \exp(-\sigma t)M\left(\sup_{0 \leq s \leq jT_2+\tau} b(s)\right)\left(\Theta_1|z_0| + \Theta_2|w_0| + \Theta_3\|x_0\|_r + \Theta_4\|u_0\|_{r+\tau}\right)$$
$$+ M\left(\sup_{0 \leq s \leq jT_2+\tau} b(s)\right)\Theta_5 \sup_{0 \leq s \leq t}(\exp(-\sigma(t-s))|\xi(s)|) + M\left(\sup_{0 \leq s \leq jT_2+\tau} b(s)\right)\Theta_6 \sup_{0 \leq s \leq t}(\exp(-\sigma(t-s))|d(s)|)$$

$$\tag{3.20}$$

where $j = \min\{j \in Z^+ : jT_2 \geq r + T_1\}$, $M(\rho) := \left(\frac{7(1+\Gamma)\exp(\beta T_2)}{\sqrt{1-\exp(-2\omega T_1 \exp(-\rho))}}\right)^{g\left(j+\frac{\tau}{T_2}\right)}$ for all $\rho \geq 0$,

$g(t) := \min\{k \in Z^+ : t \leq k\}$, $\omega := \frac{1}{2}\max\left(L(n+1)+2+2n\max_{i=1,...,n}(\theta^{2i}p_i^2), 1+L^2\right)$, $\Gamma := K\frac{((nL+1)T)^{l+1}}{1-(nL+1)T} + \exp\left(\frac{(n+1)L+3}{2}(r+\tau)\right)$

and $\beta := \omega + \frac{(n+1)L+3}{2}$.

It should be emphasized that inequality (3.19) holds for sufficiently large integers $l, m \geq 1$ and sufficiently small holding period $T_2 > 0$. The reader should notice that since $T = \frac{r+\tau}{m}$ the selection of sufficiently large integers $l, m \geq 1$ makes the term $K|k|\frac{((nL+1)T)^{l+1}}{1-(nL+1)T}$ sufficiently small: first select an integer $m \geq 1$ so that $(nL+1)(r+\tau) < m$ and then (since $K := K(m) \geq 0$ is independent of $l \geq 1$; see Proposition 3.1) we can select a sufficiently large integer $l \geq 1$ so that $K|k|\frac{((nL+1)T)^{l+1}}{1-(nL+1)T}$ becomes appropriately small.

The proof of Theorem 3.2 is based on the following technical lemmas. Their proofs are given in the Appendix.

**Lemma 3.3 (Bound on Observer State):** *For every* $(x,u,z_0,w_0) \in C^0([-r,+\infty);\Re^n) \times L^\infty_{loc}([-r-\tau,+\infty);\Re) \times \Re^n \times \Re$, $(\xi,b) \in L^\infty_{loc}(\Re_+;\Re \times \Re_+)$ *the solution* $(z(t),w(t)) \in \Re^n \times \Re$ *of the hybrid system (3.5) with initial condition* $(z(0),w(0)) = (z_0,w_0) \in \Re^n \times \Re$ *and corresponding to inputs* $(\xi,b) \in L^\infty_{loc}(\Re_+;\Re \times \Re_+)$, $(x,u) \in C^0([-r,+\infty);\Re^n) \times L^\infty_{loc}([-r-\tau,+\infty);\Re)$ *exists for all* $t \geq 0$ *and satisfies the following inequality:*

$$|z(t)|^2 + w^2(t) \leq \exp(2\omega t)\left(|z_0|^2 + |w_0|^2 + \frac{\left(\sup_{0 \leq s \leq t}|x(s-r)| + \sup_{0 \leq s \leq t}|\xi(s)|\right)^2}{1 - \exp\left(-2\omega T_1 \exp\left(-\sup_{0 \leq s \leq t} b(s)\right)\right)} + \frac{1}{2\omega}\sup_{0 \leq s < t}|u(s-r-\tau)|^2\right) \quad (3.21)$$

*where* $\omega := \frac{1}{2}\max\left(L(n+1) + 2 + 2n\max_{i=1,\dots,n}(\theta^{2i}p_i^2), 1 + L^2\right)$.

**Lemma 3.4 (Closed-Loop Solution Exists for all Times):** *For every* $(x_0,u_0,z_0,w_0) \in C^0([-r,0];\Re^n) \times L^\infty([-r-\tau,0);\Re) \times \Re^n \times \Re$, $(\xi,b,d) \in L^\infty_{loc}(\Re_+;\Re \times \Re_+ \times \Re^n)$ *the solution* $(T_r(t)x, \breve{T}_{r+\tau}(t)u, z(t), w(t)) \in C^0([-r,0];\Re^n) \times L^\infty([-r-\tau,0);\Re) \times \Re^n \times \Re$ *of the closed-loop system (1.4), (3.5) and (3.16) with initial condition* $\breve{T}_{r+\tau}(0)u = u_0 \in L^\infty([-r-\tau,0);\Re)$, $T_r(0)x = x_0 \in C^0([-r,0];\Re^n)$, $(z(0),w(0)) = (z_0,w_0) \in \Re^n \times \Re$ *and corresponding to inputs* $(\xi,b,d) \in L^\infty_{loc}(\Re_+;\Re \times \Re_+ \times \Re^n)$ *exists for all* $t \geq 0$ *and satisfies the following estimate:*

$$\sup_{0 \leq s \leq t}(|z(s)| + |w(s)|) + \sup_{-r \leq s \leq t}(|x(s)|) + \sup_{-r-\tau \leq s < t}(|u(s)|) \leq$$

$$\left(\frac{7(1+\Gamma)\exp(\beta T_2)}{\sqrt{1 - \exp\left(-2\omega T_1 \exp\left(-\sup_{0 \leq s \leq t} b(s)\right)\right)}}\right)^{g\left(\frac{t}{T_2}\right)} \left(|z_0| + |w_0| + \|x_0\|_r + \|u_0\|_{r+\tau} + \sup_{0 \leq s \leq t}|\xi(s)| + G\sup_{0 \leq s \leq t}|d(s)|\right) \quad (3.22)$$

*where* $g(t) := \min\{k \in Z^+ : t \leq k\}$, $\beta := \omega + \frac{(n+1)L+3}{2}$, $\omega := \frac{1}{2}\max\left(L(n+1) + 2 + 2n\max_{i=1,\dots,n}(\theta^{2i}p_i^2), 1 + L^2\right)$ *and*

$\Gamma := K\frac{((nL+1)T)^{l+1}}{1-(nL+1)T} + \exp\left(\frac{(n+1)L+3}{2}(r+\tau)\right)$.

**Lemma 3.5 (Convergence of Observer Estimate for Fast Sampling and High Observer Gain):** *Let* $Q \in \Re^{n \times n}$ *be a symmetric positive definite matrix that satisfies* $Q(A+pc') + (A'+cp')Q + 2qI \leq 0$ *for certain constant* $q > 0$. *Suppose that (3.17), (3.18) hold for certain constant* $a > 0$ *satisfying* $a|x|^2 \leq x'Qx$ *for all* $x \in \Re^n$. *Then there exist constants* $\sigma > 0$, $A_i > 0$ $(i=1,\dots,4)$, *which are independent of* $T_2 > 0$ *and* $l,m$, *such that for every* $(x_0,u_0,z_0,w_0) \in C^0([-r,0];\Re^n) \times L^\infty([-r-\tau,0);\Re) \times \Re^n \times \Re$, $(\xi,b,d) \in L^\infty_{loc}(\Re_+;\Re \times \Re_+ \times \Re^n)$ *the solution* $(T_r(t)x, \breve{T}_{r+\tau}(t)u, z(t), w(t)) \in C^0([-r,0];\Re^n) \times L^\infty([-r-\tau,0);\Re) \times \Re^n \times \Re$ *of the closed-loop system (1.4), (3.5) and (3.16) with initial condition* $\breve{T}_{r+\tau}(0)u = u_0 \in L^\infty([-r-\tau,0);\Re)$, $T_r(0)x = x_0 \in C^0([-r,0];\Re^n)$, $(z(0),w(0))=(z_0,w_0)\in\Re^n\times\Re$ *and corresponding to inputs* $(\xi,b,d)\in L^\infty_{loc}(\Re_+;\Re\times\Re_+\times\Re^n)$ *satisfies the following estimate for all* $t\geq r+T_1$:

$$\begin{aligned}|z(t)-x(t-r)|&\leq A_1\exp(-\sigma(t-r))|z(r)-x(0)|+A_2\sup_{0\leq s\leq t}(\exp(-\sigma(t-s))|\xi(s)|)\\&+A_3\sup_{r\leq s\leq r+T_1}(\exp(-\sigma(t-s))|w(s)-x_1(s-r)|)+A_4\sup_{0\leq s\leq t}(\exp(-\sigma(t-s))|d(s)|)\end{aligned}\quad(3.23)$$

**Lemma 3.6 (Zero-Order Hold Control Close to Nominal Control if Sampling is Fast and Approximate Predictor is Accurate):** *Let* $Q\in\Re^{n\times n}$ *be a symmetric positive definite matrix that satisfies* $Q(A+pc')+(A'+cp')Q+2qI\leq 0$ *for certain constant* $q>0$ *and certain* $p\in\Re^n$. *Suppose that inequalities (3.17), (3.18), (3.19) hold for certain constants* $a>0$ *satisfying* $a|x|^2\leq x'Qx$ *for all* $x\in\Re^n$, $0<K_1\leq K_2$ *satisfying* $K_1|x|^2\leq x'Px\leq K_2|x|^2$ *for all* $x\in\Re^n$, $K>0$ *being the constant involved in (3.13) and* $T=\frac{r+\tau}{m}$. *Define* $j=\min\{i\in Z^+:jT_2\geq r+T_1\}$. *Then for all sufficiently small* $\sigma>0$ *and for all* $(x_0,u_0,z_0,w_0)\in C^0([-r,0];\Re^n)\times L^\infty([-r-\tau,0);\Re)\times\Re^n\times\Re$, $(\xi,b,d)\in L^\infty_{loc}(\Re_+;\Re\times\Re_+\times\Re^n)$ *(independent of* $\sigma>0$*) the solution* $(T_r(t)x,\breve{T}_{r+\tau}(t)u,z(t),w(t))\in C^0([-r,0];\Re^n)\times L^\infty([-r-\tau,0);\Re)\times\Re^n\times\Re$ *of the closed-loop system (1.4), (3.5) and (3.16) with initial condition* $\breve{T}_{r+\tau}(0)u=u_0\in L^\infty([-r-\tau,0);\Re)$, $T_r(0)x=x_0\in C^0([-r,0];\Re^n)$, $(z(0),w(0))=(z_0,w_0)\in\Re^n\times\Re$ *and corresponding to inputs* $(\xi,b,d)\in L^\infty_{loc}(\Re_+;\Re\times\Re_+\times\Re^n)$ *satisfies the following estimate for all* $t\geq jT_2+\tau$:

$$\begin{aligned}&\eta\exp(\sigma t)|u(t-\tau)-k'x(t)|\\&\leq C|k|\exp(\sigma(T_2+r+\tau))\sup_{jT_2-r\leq s<jT_2+\tau}(\exp(\sigma s)|u(s-\tau)-k'x(s)|)\\&+A_1\exp(\sigma(T_2+r+\tau))|k|(C+\exp((nL+1)(r+\tau)))|z(r)-x(0)|\\&+A_2\exp(\sigma(T_2+\tau))|k|(C+\exp((nL+1)(r+\tau)))\sup_{0\leq s\leq t}(\exp(\sigma s)|\xi(s)|)\\&+A_3\exp(\sigma(T_2+\tau))|k|(+\exp((nL+1)(r+\tau)))\sup_{r\leq s\leq r+T_1}(\exp(\sigma s)|w(s)-x_1(s-r)|)\\&+|k|\exp(\sigma(T_2+\tau))[A_4C+\exp((nL+1)(r+\tau))(A_4+G(r+\tau)\exp(\sigma r))+T_2G\exp(-\sigma\tau)]\sup_{0\leq s\leq t}(\exp(\sigma s)|d(s)|)\\&+|k|\exp(\sigma T_2)[C(1+|k|)\exp(\sigma(r+\tau))+T_2(nL+1+|k|)]\sup_{-r\leq s\leq t}(\exp(\sigma s)|x(s)|)\end{aligned}\quad(3.24)$$

*where* $\eta:=1-|k|T_2-C|k|\exp(\sigma(T_2+r+\tau))$, $C:=K\frac{((nL+1)T)^{l+1}}{1-(nL+1)T}$ *and* $A_i>0$ ($i=1,\ldots,4$) *are the constants involved in (3.23).*

We now provide the proof of Theorem 3.2.

**Proof of Theorem 3.2:** Let $\sigma>0$ be sufficiently small such that $\sqrt{\frac{b'Pb}{2K_1}}|k|\exp(\sigma T_2)[C(1+|k|)\exp(\sigma(r+\tau))+T_2(nL+1+|k|)]<\eta\mu$, $\sigma\leq\frac{\mu}{2}$ and such that inequalities (3.23), (3.24) hold. The existence of sufficiently small $\sigma>0$ satisfying $\sqrt{\frac{b'Pb}{2K_1}}|k|\exp(\sigma T_2)[C(1+|k|)\exp(\sigma(r+\tau))+T_2(nL+1+|k|)]<\eta\mu$ is a direct consequence of (3.19). Define $V(t)=x'(t)Px(t)$. Using (3.15) we obtain the following differential inequality for almost all $t\geq 0$:

$$\dot{V}(t) \leq -2\mu V(t) + \frac{b'Pb}{2\mu}|u(t-\tau) - k'x(t)|^2 + 2\gamma|d(t)|^2 \tag{3.25}$$

The above differential inequality directly gives the following estimate for all $t > 0$:

$$|x(t)| \leq \sqrt{\frac{K_2}{K_1}}\exp(-\mu t)|x(0)| + \frac{1}{\mu}\sqrt{\frac{b'Pb}{2K_1}}\sup_{0 \leq s < t}\left(\exp\left(-\frac{\mu}{2}(t-s)\right)|u(s-\tau) - k'x(s)|\right)$$
$$+ \sqrt{\frac{2\gamma}{\mu K_1}}\sup_{0 \leq s \leq t}\left(\exp\left(-\frac{\mu}{2}(t-s)\right)|d(s)|\right) \tag{3.26}$$

where $0 < K_1 \leq K_2$ are constants satisfying $K_1|x|^2 \leq x'Px \leq K_2|x|^2$ for all $x \in \Re^n$. Using the inequality $\sigma \leq \frac{\mu}{2}$ we conclude that the following inequality holds for all $t > 0$:

$$|x(t)| \leq \sqrt{\frac{K_2}{K_1}}\exp(-\sigma t)|x(0)| + \frac{1}{\mu}\sqrt{\frac{b'Pb}{2K_1}}\sup_{0 \leq s < t}(\exp(-\sigma(t-s))|u(s-\tau) - k'x(s)|) + \sqrt{\frac{2\gamma}{\mu K_1}}\sup_{0 \leq s < t}(\exp(-\sigma(t-s))|d(s)|) \tag{3.27}$$

Notice that inequality (3.27) implies the following inequality for all $t > 0$:

$$\sup_{0 \leq s \leq t}(\exp(\sigma s)|x(s)|) \leq \sqrt{\frac{K_2}{K_1}}|x(0)| + \frac{1}{\mu}\sqrt{\frac{b'Pb}{2K_1}}\sup_{0 \leq s < t}(\exp(\sigma s)|u(s-\tau) - k'x(s)|) + \sqrt{\frac{2\gamma}{\mu K_1}}\sup_{0 \leq s < t}(\exp(\sigma s)|d(s)|) \tag{3.28}$$

Combining (3.28) and (3.24), we obtain the following inequality for all $t \geq jT_2 + \tau$, where $j = \min\{j \in Z^+ : jT_2 \geq r + T_1\}$:

$$\sup_{0 \leq s \leq t}(\exp(\sigma s)|x(s)|) \leq \sqrt{\frac{K_2}{K_1}}|x(0)| + \sqrt{\frac{2\gamma}{\mu K_1}}\sup_{0 \leq s \leq t}(\exp(\sigma s)|d(s)|)$$
$$+ \frac{1}{\mu}\sqrt{\frac{b'Pb}{2K_1}}(1 + \eta^{-1}C|k|\exp(\sigma(T_2 + r + \tau)))\sup_{0 \leq s < jT_2 + \tau}(\exp(\sigma s)|u(s-\tau) - k'x(s)|)$$
$$+ \frac{1}{\eta\mu}\sqrt{\frac{b'Pb}{2K_1}}A_1\exp(\sigma(T_2 + r + \tau))|k|(C + \exp((nL+1)(r+\tau)))|z(r) - x(0)|$$
$$+ \frac{1}{\eta\mu}\sqrt{\frac{b'Pb}{2K_1}}A_2\exp(\sigma(T_2 + \tau))|k|(C + \exp((nL+1)(r+\tau)))\sup_{0 \leq s \leq t}(\exp(\sigma s)|\xi(s)|)$$
$$+ \frac{1}{\eta\mu}\sqrt{\frac{b'Pb}{2K_1}}A_3\exp(\sigma(T_2 + \tau))|k|(C + \exp((nL+1)(r+\tau)))\sup_{r \leq s \leq r+T_1}(\exp(\sigma s)|w(s) - x_1(s-r)|)$$
$$+ \frac{1}{\eta\mu}\sqrt{\frac{b'Pb}{2K_1}}|k|\exp(\sigma(T_2 + \tau))[A_4C + \exp((nL+1)(r+\tau))(A_4 + G(r+\tau)\exp(\sigma r)) + T_2G\exp(-\sigma\tau)]\sup_{0 \leq s \leq t}(\exp(\sigma s)|d(s)|)$$
$$+ \frac{1}{\eta\mu}\sqrt{\frac{b'Pb}{2K_1}}|k|\exp(\sigma T_2)[C(1+|k|)\exp(\sigma(r+\tau)) + T_2(nL+1+|k|)]\|x_0\|_r$$
$$+ \frac{1}{\eta\mu}\sqrt{\frac{b'Pb}{2K_1}}|k|\exp(\sigma T_2)[C(1+|k|)\exp(\sigma(r+\tau)) + T_2(nL+1+|k|)]\sup_{0 \leq s \leq t}(\exp(\sigma s)|x(s)|)$$

It is clear from the above inequality that there exist constants $B_i > 0$ ($i = 1,...,6$) so that the following inequality holds for all $t \geq jT_2 + \tau$, where $j = \min\{j \in Z^+ : jT_2 \geq r + T_1\}$:

$$\sup_{0 \leq s \leq t}(\exp(\sigma s)|x(s)|) \leq B_1\|x_0\|_r + B_2\sup_{0 \leq s \leq t}(\exp(\sigma s)|d(s)|) + B_3\sup_{0 \leq s < jT_2 + \tau}(\exp(\sigma s)|u(s-\tau) - k'x(s)|)$$
$$+ B_4|z(r) - x(0)| + B_5\sup_{0 \leq s \leq t}(\exp(\sigma s)|\xi(s)|) + B_6\sup_{r \leq s \leq r+T_1}(\exp(\sigma s)|w(s) - x_1(s-r)|) \tag{3.28}$$

provided that $\sqrt{\frac{b'Pb}{2K_1}}|k|\exp(\sigma T_2)[C(1+|k|)\exp(\sigma(r+\tau))+T_2(nL+1+|k|)] < \eta\mu$, where

$\eta := 1-|k|T_2 - C|k|\exp(\sigma(T_2+r+\tau))$, $C := K\frac{((nL+1)T)^{l+1}}{1-(nL+1)T}$. Combining inequalities (3.24), (3.28), (3.23), (3.22) and inequality (A.18) in the proof of Lemma 3.5, we obtain the existence of constants $\Theta_i > 0$ ($i = 1,...,6$) satisfying inequality (3.20). The proof is complete. ◁

**Remark 3.7:**

(a) Small errors in the measurement delay $r \geq 0$ can be handled. In order to see this, notice that a small error in the measurement delay $r \geq 0$ induces a measurement error $\xi(t) = x(t-\hat{r}) - x(t-r)$, where $r \geq 0$ is the assumed value of the measurement delay and $\hat{r} \geq 0$ is the actual value of the measurement delay. If $|r-\hat{r}|$ is sufficiently small then the measurement error $\xi(t) = x(t-\hat{r}) - x(t-r)$ can be rendered sufficiently small so that exponential convergence is preserved. More specifically, there exist constants $\Delta_i > 0$ ($i = 1,...,4$) such that for every $(x_0, u_0, z_0, w_0) \in C^0([-R,0];\Re^n) \times L^\infty([-r-\tau,0);\Re) \times \Re^n \times \Re$, $(\xi, b, d) \in L^\infty_{loc}(\Re_+;\Re \times \Re_+ \times \Re^n)$, where $R := \max(r, \hat{r})$, whenever the solution $(T_r(t)x, \breve{T}_{r+\tau}(t)u, z(t), w(t)) \in C^0([-R,0];\Re^n) \times L^\infty([-r-\tau,0);\Re) \times \Re^n \times \Re$ of the closed-loop system (1.4), (3.5) and (3.16) with initial condition $\breve{T}_{r+\tau}(0)u = u_0 \in L^\infty([-r-\tau,0);\Re)$, $T_R(0)x = x_0 \in C^0([-R,0];\Re^n)$, $(z(0), w(0)) = (z_0, w_0) \in \Re^n \times \Re$ and corresponding to inputs $(\xi, b, d) \in L^\infty_{loc}(\Re_+;\Re \times \Re_+ \times \Re^n)$ exists, the solution satisfies the following estimate for all $t \geq 0$:

$$\exp(\sigma t)|\xi(t)| = \exp(\sigma t)|x(t-r) - x(t-\hat{r})|$$
$$\leq \Delta_1 \|x_0\|_R + |r-\hat{r}|\Delta_2 \sup_{0 \leq s \leq t}(\exp(\sigma s)|d(s)|) + |r-\hat{r}|\Delta_3 \sup_{-R \leq s \leq t}(\exp(\sigma s)|x(s)|) \quad (3.29)$$
$$+ |r-\hat{r}|\Delta_4 \sup_{-\tau \leq s < t}(\exp(\sigma s)|u(s)|)$$

Combining (3.29) with (3.20) we conclude that the following estimate for all $t \geq 0$ for which the solution exists:

$$\frac{\sup_{0 \leq s \leq t}(\exp(\sigma s)(|z(s)|+|w(s)|+\|T_r(s)x\|_r)) + \sup_{0 \leq s \leq t}(\exp(\sigma s)\|\breve{T}_{r+\tau}(s)u\|_{r+\tau})}{M\left(\sup_{0 \leq s \leq jT_2+\tau} b(s)\right)} \leq$$

$$2\Theta_1|z_0| + 2\Theta_2|w_0| + 2(\Theta_3 + \Theta_5\Delta_1)\|x_0\|_R + 2\Theta_4\|u_0\|_{r+\tau}$$
$$+ 2|r-\hat{r}|\Delta_3\Theta_5 \sup_{0 \leq s \leq t}(\exp(\sigma s)|x(s)|) + 2|r-\hat{r}|\Delta_4\Theta_5 \sup_{-\tau \leq s < t}(\exp(\sigma s)|u(s)|) + 2(\Theta_6 + |r-\hat{r}|\Delta_2\Theta_5) \sup_{0 \leq s \leq t}(\exp(\sigma s)|d(s)|)$$

From the above inequality, the solution of the closed-loop system (1.4), (3.5) and (3.16) exists for all $t \geq 0$ provided that the inequality $2|r-\hat{r}|\Theta_5 \max(\Delta_3, \Delta_4) M\left(\sup_{0 \leq s \leq jT_2+\tau} b(s)\right) < 1$ holds. Moreover, the solution of the closed-loop system (1.4), (3.5) and (3.16) converges exponentially to zero.

(b) For the case that the input can be continuously adjusted, a similar result to Theorem 3.2 can be proved. The controller will be a combination of the prediction-based controller proposed in [20], the sampled-data observer (3.5) and the approximate predictor mapping defined by (3.12).

# 4. Specialization to Linear Time Invariant Systems

For the LTI case (1.6), where the pair of matrices $A \in \Re^{n \times n}$, $B \in \Re^n$ is stabilizable and the output is given by

$$y(\tau_i) = c'x(\tau_i - r) + \xi(\tau_i), \quad i \in Z^+ \tag{4.1}$$

where $\{\tau_i\}_{i=0}^\infty$ is a partition of $\Re^+$ with $\sup_{i \geq 0}(\tau_{i+1} - \tau_i) \leq T_1$ and the pair of matrices $A \in \Re^{n \times n}$, $c \in \Re^n$ is a detectable pair, we apply the observer-based predictor stabilization scheme described in Section 3. There exist vectors $k \in \Re^n$ and $p \in \Re^n$ such that the matrices $A + Bk'$ and $A + pc'$ are Hurwitz matrices. Moreover, the predictor mapping that relates $x(t - r)$ with $x(t + \tau)$ is given by the explicit expression

$$\Phi(x, u) := \exp(A(r + \tau))x + \int_{-r-\tau}^{0} \exp(-As) Bu(s) ds$$

Notice that the above prediction scheme is exact (not approximate) for the case $d \equiv 0$. Therefore, we prove the following corollary in exactly the same way as in Theorem 3.2.

**Corollary 4.1:** *Assume that there exist vectors $k \in \Re^n$, $p \in \Re^n$ such that the matrices $A + Bk'$, $A + pc'$ are Hurwitz matrices. For sufficiently small holding period $T_2 > 0$ and for sufficiently small sampling period $T_1 > 0$, there exist constants $\Theta_i > 0$ ($i = 1, ..., 7$) and $\sigma, \omega, \beta > 0$ such that for every $(x_0, u_0, z_0, w_0) \in C^0([-r, 0]; \Re^n) \times L^\infty([-r - \tau, 0); \Re) \times \Re^n \times \Re$, $(\xi, b, d) \in L^\infty_{loc}(\Re_+; \Re \times \Re_+ \times \Re^n)$ the solution $(T_r(t)x, \breve{T}_{r+\tau}(t)u, z(t), w(t)) \in C^0([-r, 0]; \Re^n) \times L^\infty([-r - \tau, 0); \Re) \times \Re^n \times \Re$ of the closed-loop system consisting of (1.6) with*

$$\begin{aligned} \dot{z}(t) &= Az(t) + Bu(t - r - \tau) + p(c'z(t) - w(t)) \\ \dot{w}(t) &= c'Az(t) + c'Bu(t - r - \tau), \quad t \in [\tau_i, \tau_{i+1}) \\ w(\tau_{i+1}) &= y(\tau_{i+1}) = c'x(\tau_{i+1} - r) + \xi(\tau_{i+1}) \\ \tau_{i+1} &= \tau_i + T_1 \exp(-b(\tau_i)), \tau_0 = 0 \end{aligned} \tag{4.2}$$

$$u(t) = k' \exp(A(r + \tau)) z(iT_2) + \int_{-r-\tau}^{0} k' \exp(-As) Bu(iT_2 + s) ds, \text{ for } t \in [iT_2, (i+1)T_2) \tag{4.3}$$

*and initial condition $\breve{T}_{r+\tau}(0)u = u_0 \in L^\infty([-r - \tau, 0); \Re)$, $T_r(0)x = x_0 \in C^0([-r, 0]; \Re^n)$, $(z(0), w(0)) = (z_0, w_0) \in \Re^n \times \Re$ and corresponding to inputs $(\xi, b, d) \in L^\infty_{loc}(\Re_+; \Re \times \Re_+ \times \Re^n)$ satisfies inequality (3.20) for all $t \geq 0$, where $j = \min\{j \in Z^+ : jT_2 \geq r + T_1\}$,*

$$M(\rho) := \left( \frac{\Theta_7 \exp(\beta T_2)}{\sqrt{1 - \exp(-\omega T_1 \exp(-\rho))}} \right)^{g\left(j + \frac{\tau}{T_2}\right)} \text{ for all } \rho \geq 0 \text{ and } g(t) := \min\{k \in Z^+ : t \leq k\}.$$

The advantage of the sampled-data feedback stabilizer (4.2), (4.3) compared to other sampled-data stabilizers for (1.6) (see for example [29]) is that the closed-loop system (1.6), (4.2), (4.3) is completely insensitive to perturbations of the sampling schedule (this is guaranteed by inequality (3.20) and the fact that possible perturbations of the sampling schedule are quantified by means of the input $b \in L^\infty_{loc}(\Re_+; \Re_+)$).

# 5. Illustrative Example

In this section we consider the following two dimensional system

$$\dot{x}_1(t) = f(x_1(t)) + x_2(t), \; \dot{x}_2(t) = u(t - \tau) \tag{5.1}$$

where $f(x) = \text{sgn}(x)\dfrac{x^2}{\sqrt{1+x^2}}$. For this function we have $f'(x) = \dfrac{|x|(2+x^2)}{(1+x^2)\sqrt{1+x^2}}$, $\sup_{x \in \Re}|f'(x)| = \dfrac{4\sqrt{2}}{3\sqrt{3}} \approx 1.088662$ and consequently system (5.1) is of the form (1.4) and satisfies the global Lipschitz assumption made in Section 3. The one-dimensional version of system (5.1) was studied in [18], where it was shown that a nonlinear predictor scheme was necessary for its stabilization. Here, we study system (5.1) with output available at discrete time instants:

$$y(t) = x_1(iT_1 - r), \text{ for } t \in [iT_1, (i+1)T_1), \; i \in Z^+ \tag{5.2}$$

where $T_1 > 0$ is the sampling period and $r \geq 0$ is the measurement delay. The input $u(t)$ is applied with zero-order hold with holding period $T_2 > 0$. Theorem 3.2 implies that there exist constants $\Theta_i > 0$ ($i=1,...,4$) and $\sigma > 0$ such that for every $(x_0, u_0, z_0, w_0) \in C^0([-r,0];\Re^n) \times L^\infty([-r-\tau,0);\Re) \times \Re^n \times \Re$ the solution $(T_r(t)x, \breve{T}_{r+\tau}(t)u, z(t), w(t)) \in C^0([-r,0];\Re^n) \times L^\infty([-r-\tau,0);\Re) \times \Re^n \times \Re$ of the closed-loop system (5.1) with

$$\begin{aligned}
\dot{z}_1(t) &= f(z_1(t)) + z_2(t) - 3\theta(z_1(t) - w(t)) \\
\dot{z}_2(t) &= -3\theta^2(z_1(t) - w(t)) + u(t - r - \tau) \\
\dot{w}(t) &= f(z_1(t)) + z_2(t) \quad , \quad t \in [iT_1, (i+1)T_1), i \in Z^+ \\
w((i+1)T_1) &= y((i+1)T_1) = x_1((i+1)T_1 - r)
\end{aligned} \tag{5.3}$$

$$u(t) = k'\Phi_{l,m}(z(iT_2), \breve{T}_{r+\tau}(iT_2)u), \text{ for } t \in [iT_2, (i+1)T_2) \tag{5.4}$$

where $l, m \geq 1$ are integers, the operator $\Phi_{l,m} : \Re^2 \times L^\infty([-r-\tau,0);\Re) \to \Re^2$ is defined by (3.12), $k = -(15,9)' \in \Re^2$ and initial condition $\breve{T}_{r+\tau}(0)u = u_0 \in L^\infty([-r-\tau,0);\Re)$, $T_r(0)x = x_0 \in C^0([-r,0];\Re^n)$, $(z(0), w(0)) = (z_0, w_0) \in \Re^n \times \Re$ satisfies the following inequality for all $t \geq 0$:

$$|z(t)| + |w(t)| + \|T_r(t)x\|_r + \|\breve{T}_{r+\tau}(t)u\|_{r+\tau} \leq \exp(-\sigma t)(\Theta_1|z_0| + \Theta_2|w_0| + \Theta_3\|x_0\|_r + \Theta_4\|u_0\|_r) \tag{5.5}$$

provided that $l, m$ are sufficiently large positive integers, $\theta \geq 1$ is sufficiently large and the sampling period $T_1 > 0$ and holding period $T_2 > 0$ are sufficiently small.

The closed-loop system (5.1), (5.3), (5.4) was tested numerically for $r = \tau = 1/4$. It was found that the selection

$$l = m = 1, \; \Phi_{l,m}(z_1, z_2, u) = \frac{1}{2}\begin{bmatrix} 2z_1 + z_2 + f(z_1) \\ 2z_2 + 2\int_{-1/2}^{0} u(s)ds \end{bmatrix}, \; \theta = 1, \; T_2 = \frac{1}{100}, \; T_1 = 3T_2 = \frac{3}{100} \tag{5.6}$$

was appropriate in order to guarantee exponential stability for the closed-loop system. Figures 1 and 2 show the time evolution of the state and the input for initial conditions $x_1(s) = x_2(s) = 1$ for $s \in \left[-\dfrac{1}{4}, 0\right]$, $u(s) = -2$ for $s \in \left[-\dfrac{1}{2}, 0\right)$ and $z_1(0) = z_2(0) = w(0) = 0$. It is clearly shown that all variables converge exponentially to zero.

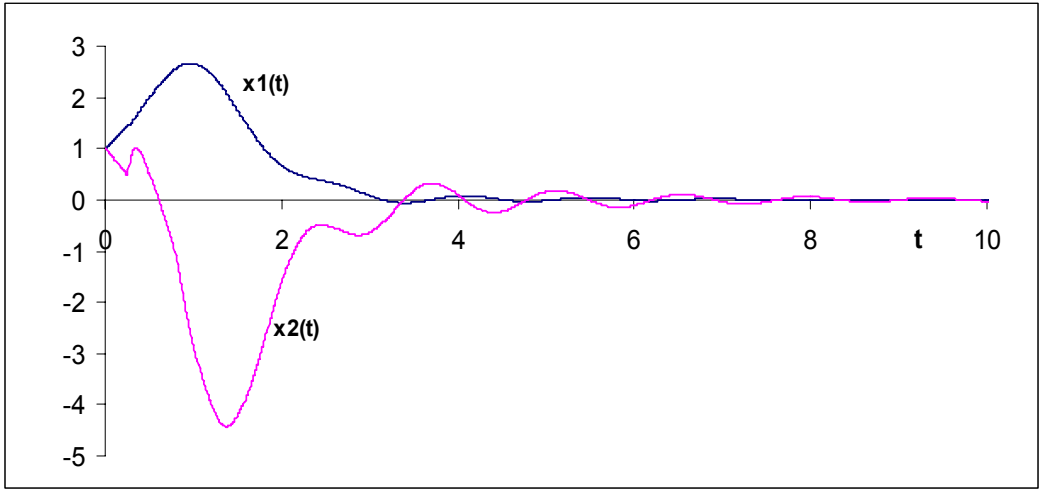

**Figure 1:** Time evolution of the state $(x_1(t), x_2(t))$ of the closed-loop system (5.1), (5.3), (5.4), (5.6) with initial conditions $x_1(s) = x_2(s) = 1$ for $s \in \left[-\frac{1}{4}, 0\right]$, $u(s) = -2$ for $s \in \left[-\frac{1}{2}, 0\right)$ and $z_1(0) = z_2(0) = w(0) = 0$

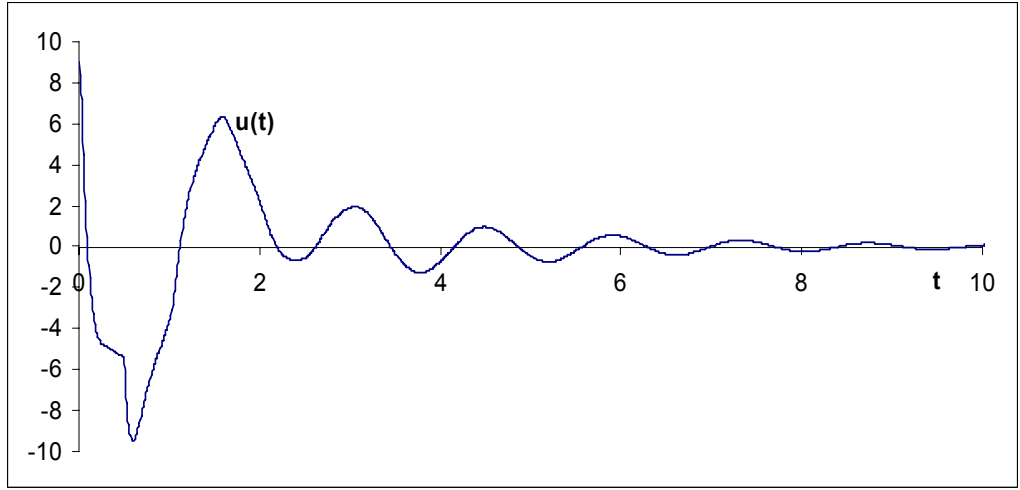

**Figure 2:** Time evolution of the input $u(t)$ for the closed-loop system (5.1), (5.3), (5.4), (5.6) with initial conditions $x_1(s) = x_2(s) = 1$ for $s \in \left[-\frac{1}{4}, 0\right]$, $u(s) = -2$ for $s \in \left[-\frac{1}{2}, 0\right)$ and $z_1(0) = z_2(0) = w(0) = 0$

## 6. Concluding Remarks

We have expanded the applicability of delay-compensating stabilizing feedback to nonlinear systems where only output measurement is available and where such measurement is subject to long delays. Our designs employ either exact or approximate predictor maps. We perform state estimation using either reconstruction maps that generate the state in a finite number of steps from output and input data, or using high-gain sampled-data observers. Our results are global, and guarantee input-to-state stability in the presence of disturbances for globally Lipschitz systems, provided the sampling/holding periods are sufficiently short. Numerous relevant open problems remain that include multiple delays on inputs, states, and in the output map or quantization issues (as in [4,5,6]), or the possible use of emulation-based observers (as in [2]).

# Appendix

**Proof of Theorem 2.3:** By virtue of the Fact, inequality (2.3), the fact that the reconstruction mapping $R: U^{p+1} \times \Re^{kp} \times \Re^k \to \Re^n$ is a continuous function with $R(0,0) = 0$, in conjunction with (2.24), (2.25), there exists a function $b \in K_\infty$ such that

$$|u_i| \leq b\left(\max_{j=i-p,\ldots,i} |y_j| + \max_{j=i-p-l-1,\ldots,i-1} |u_j|\right), \text{ for all } i \in Z^+ \quad \text{(A.1)}$$

Let $(x_0, u_0) \in C^0([-r-pT,0]; \Re^n) \times L^\infty([-(p+l+1)T,0); U)$ be arbitrary and consider the solution $(x(t), u(t)) \in \Re^n \times U$ of the closed-loop system (1.1), (1.2), (1.3), (2.24), (2.25) with $T_1 = T_2 = T$, initial condition $\breve{T}_{r+\tau}(0)u = u_0 \in L^\infty([-(p+l+1)T,0); U)$, $T_r(0)x = x_0 \in C^0([-r-pT,0]; \Re^n)$. Using (A.1) and the Fact, we can show that the solution of the closed-loop system (1.1), (1.2), (1.3) with $T_1 = T_2 = T$, (2.24), (2.25) exists for all $t \geq 0$. Moreover, using (A.1) and the Fact we can construct inductively a function $\tilde{b} \in K_\infty$ such that

$$\|T_{r+pT}(t)x\|_{r+pT} + \|\breve{T}_{(p+l+1)T}(t)u\|_{(p+l+1)T} \leq \tilde{b}\left(\|x_0\|_{r+pT} + \|u_0\|_{(p+l+1)T}\right), \text{ for all } t \in [(p+l+1)T+\tau] \quad \text{(A.2)}$$

Notice that for every $i \in Z^+$ with $i \geq p+l+1$, (2.20) holds and $u_i = k(x(iT+\tau))$. Hypothesis (H2) and (2.3) imply that there exists $\sigma \in KL$ such that

$$|x(t)| + |u(t)| \leq \sigma\left(|x((p+l+1)T+\tau)|, t-(p+l+1)T-\tau\right), \text{ for all } t \geq (p+l+1)T+\tau \quad \text{(A.3)}$$

Combining (A.2) and (A.3) we can guarantee that there exists a function $\tilde{\sigma} \in KL$ such that (2.26) holds. Finally, if the closed-loop system (2.1), (2.2) satisfies the dead-beat property of order $jT$, where $j \in Z^+$ is positive, then (2.20) implies that $x(t) = 0$ for all $t \geq (j+p+l+1)T+\tau$. The proof is complete. ◁

**Proof of Lemma 3.3:** Local existence and uniqueness follows from [19] (pages 23-27). Moreover, the analysis in [19] (pages 23-27) shows that the solution exists as long as it is bounded. In order to show that the solution remains bounded for all finite times, we consider the function $R(t) = \frac{1}{2}|z(t)|^2 + \frac{1}{2}w^2(t)$. Using algebraic manipulations and (3.2), it follows that the following differential inequality holds for almost all $t \in [\tau_i, \tau_{i+1})$ and $i \in Z^+$ for which the solution exists:

$$\dot{R}(t) \leq 2\omega R(t) + \frac{1}{2}u^2(t-r-\tau) \quad \text{(A.4)}$$

where $\omega := \frac{1}{2}\max\left(L(n+1) + 2 + 2n \max_{i=1,\ldots,n}(\theta^{2i} p_i^2), 1+L^2\right)$. Integrating the differential inequality (A.4) we obtain for all $t \in [\tau_i, \tau_{i+1})$ and $i \in Z^+$ for which the solution exists:

$$|z(t)|^2 + w^2(t) \leq \exp(2\omega(t-\tau_i))\left(|z(\tau_i)|^2 + |w(\tau_i)|^2 + \sup_{\tau_i \leq s < t}|u(s-r-\tau)|^2 \int_0^{t-\tau_i} \exp(-2\omega s)ds\right) \quad \text{(A.5)}$$

Consequently, using a standard contradiction argument and (A.5), we are able to show that:

"if for some $i \in Z^+$ the solution exists at $t = \tau_i$ then the solution exists at $t = \tau_{i+1}$"

Using induction, (A.5) and the fact that $|w(\tau_i)| \leq \sup_{0 \leq s \leq \tau_i} |x(s-r)| + \sup_{0 \leq s \leq \tau_i} |\xi(s)|$ for all $i \in Z^+$ with $i \geq 1$, we show that the following inequality holds for all $i \in Z^+$ with $i \geq 2$:

$$\exp(-2\omega\tau_i)|z(\tau_i)|^2 \leq |z(0)|^2 + |w(0)|^2 + \left(\sup_{0 \leq s \leq \tau_i}|x(s-r)| + \sup_{0 \leq s \leq \tau_i}|\xi(s)|\right)^2 \sum_{k=1}^{i-1} \exp(-2\omega\tau_k) + \sup_{0 \leq s < \tau_i}|u(s-r-\tau)|^2 \int_0^{\tau_i} \exp(-2\omega s) ds \quad (A.6)$$

Inequalities (A.5), (A.6) and the fact that $|w(\tau_i)| \leq \sup_{0 \leq s \leq \tau_i}|x(s-r)| + \sup_{0 \leq s \leq \tau_i}|\xi(s)|$ for all $i \in Z^+$ with $i \geq 1$, we show that the following inequality holds for all $i \in Z^+$ and $t \in [\tau_i, \tau_{i+1})$:

$$|z(t)|^2 + w^2(t) \leq \exp(2\omega t)\left(|z(0)|^2 + |w(0)|^2 + \left(\sup_{0 \leq s \leq t}|x(s-r)| + \sup_{0 \leq s \leq t}|\xi(s)|\right)^2 \sum_{k=0}^{i} \exp(-2\omega\tau_k) + \frac{1}{2\omega}\sup_{0 \leq s < t}|u(s-r-\tau)|^2\right)$$
(A.7)

Inequality (3.21) is a direct consequence of (A.7) and the fact that $\tau_{i+1} \geq \tau_i + T_1 \exp\left(-\sup_{0 \leq s \leq t} b(s)\right)$, which holds for all $i \in Z^+$ with $t \geq \tau_i$. The proof of the claim is complete. ◁

**Proof of Lemma 3.4:** We prove the lemma by induction. More specifically, we prove the following claim for all $i \in Z^+$:

(Claim) For every $(x_0, u_0, z_0, w_0) \in C^0([-r, 0]; \Re^n) \times L^\infty([-r-\tau, 0); \Re) \times \Re^n \times \Re$, $(\xi, b, d) \in L^\infty_{loc}(\Re_+; \Re \times \Re_+ \times \Re^n)$ the solution $(T_r(t)x, \breve{T}_{r+\tau}(t)u, z(t), w(t)) \in C^0([-r, 0]; \Re^n) \times L^\infty([-r-\tau, 0); \Re) \times \Re^n \times \Re$ of the closed-loop system (1.4), (3.5) and (3.16) with initial condition $\breve{T}_{r+\tau}(0)u = u_0 \in L^\infty([-r-\tau, 0); \Re)$, $T_r(0)x = x_0 \in C^0([-r, 0]; \Re^n)$, $(z(0), w(0)) = (z_0, w_0) \in \Re^n \times \Re$ and corresponding to inputs $(\xi, b, d) \in L^\infty_{loc}(\Re_+; \Re \times \Re_+ \times \Re^n)$ exists for all $t \in [0, iT_2]$ and satisfies (3.22) for all $t \in [0, iT_2]$.

It is clear that the claim holds for $i = 0$. Next assume that the claim holds for some $i \in Z^+$. Define

$$A_i := \left(\frac{7(1+\Gamma)\exp(\beta T_2)}{\sqrt{1 - \exp\left(-2\omega T_1 \exp\left(-\sup_{0 \leq s \leq t} b(s)\right)\right)}}\right)^i \left(|z_0| + |w_0| + \|x_0\|_r + \|u_0\|_{r+\tau} + \sup_{0 \leq s \leq t}|\xi(s)| + G \sup_{0 \leq s \leq t}|d(s)|\right).$$ Using (3.14),

(3.16) and (3.22) for $t \in [0, iT_2]$, it is clear that $u(t)$ is well-defined on $[iT_2, (i+1)T_2)$ and satisfies the following inequality for all $t \in [iT_2, (i+1)T_2]$:

$$\sup_{-r-\tau \leq s < t}(|u(s)|) \leq \Gamma A_i \quad (A.8)$$

Using (3.4), (3.22) for $t \in [0, iT_2]$ and (A.8), it is clear that $x(t)$ is well-defined on $[iT_2, (i+1)T_2]$ and satisfies the following inequality for all $t \in [iT_2, (i+1)T_2]$:

$$|x(t)| \leq \left((1+\Gamma)A_i + G\sup_{0 \leq s \leq t}|d(s)|\right)\exp\left(\frac{(n+1)L+3}{2}T_2\right) \quad (A.9)$$

Using Lemma 3.3, (A.8) and (A.9), it is clear that $(z(t), w(t))$ is well-defined on $[iT_2, (i+1)T_2]$ and satisfies the following inequality for all $t \in [iT_2, (i+1)T_2]$:

$$|z(t)| + |w(t)| \leq 2\exp(\omega T_2)\exp\left(\frac{(n+1)L+3}{2}T_2\right)\frac{G\sup_{0\leq s\leq t}|d(s)| + \sup_{0\leq s\leq t}|\xi(s)| + 2(1+\Gamma)A_i}{\sqrt{1-\exp\left(-2\omega T_1\exp\left(-\sup_{0\leq s\leq t}b(s)\right)\right)}} \quad \text{(A.10)}$$

Therefore, using the definition $\beta := \omega + \frac{(n+1)L+3}{2}$ and (A.8), (A.9), (A.10), we conclude that the following inequality holds for all $t \in [iT_2, (i+1)T_2]$:

$$\sup_{0\leq s\leq t}(|z(s)|+|w(s)|) + \sup_{-r\leq s\leq t}(|x(s)|) + \sup_{-r-\tau\leq s<t}(|u(s)|) \leq \frac{7(1+\Gamma)\exp(\beta T_2)}{\sqrt{1-\exp\left(-2\omega T_1\exp\left(-\sup_{0\leq s\leq t}b(s)\right)\right)}} A_i \quad \text{(A.11)}$$

The fact that the claim holds for all $t \in [0, (i+1)T_2]$ is a direct consequence of (A.11). The proof is complete. ◁

**Proof of Lemma 3.5:** Define the quadratic error Lyapunov function $V(e) := e'\Delta_\theta^{-1}Q\Delta_\theta^{-1}e$, where $e(t) := z(t) - x(t-r)$, $\Delta_\theta := diag(\theta, \theta^2, ..., \theta^n)$. Using (3.2), (3.3), the identities $\Delta_\theta^{-1}A = \theta A\Delta_\theta^{-1}$, $c' = \theta c'\Delta_\theta^{-1}$ and the inequalities $\theta^{-i}|f_i(x_1+e_1,...,x_i+e_i) - f_i(x_1,...,x_i)| \leq L|\Delta_\theta^{-1}e|$ for $i=1,...,n$ and all $(x,e) \in \Re^n \times \Re^n$ (which follow from (3.2)), we get for $\theta \geq \max\left(1, \frac{2|Q|L\sqrt{n}}{q}\right)$ and for all $t \geq r$:

$$\dot{e}(t) = (A + \Delta_\theta pc')e(t) + \Delta_\theta p\eta(t) + \tilde{p}(x(t-r), e(t)) - diag(g_1(x(t-r),u(t-r)),...,g_n(x(t-r),u(t-r)))d(t-r) \quad \text{(A.12)}$$

$$\dot{V}(t) \leq -2\theta q|\Delta_\theta^{-1}e(t)|^2 + 2|\Delta_\theta^{-1}e(t)||Q||\Delta_\theta^{-1}\tilde{p}(x(t-r),e(t))| + 2|e'(t)\Delta_\theta^{-1}||Qp||\eta(t)| + 2\theta^{-1}G|e'(t)\Delta_\theta^{-1}||Q||d(t-r)| \leq$$

$$\leq -2\theta q|\Delta_\theta^{-1}e(t)|^2 + 2|\Delta_\theta^{-1}e(t)|^2|Q|L\sqrt{n} + \frac{\theta q}{2}|e'(t)\Delta_\theta^{-1}|^2 + \frac{4}{\theta q}|Qp|^2|\eta(t)|^2 + \frac{4G^2}{\theta^3 q}|Q|^2|d(t-r)|^2$$

$$\leq -\frac{\theta q}{2}|\Delta_\theta^{-1}e(t)|^2 + \frac{4}{\theta q}|Qp|^2|\eta(t)|^2 + \frac{4G^2}{\theta^3 q}|Q|^2|d(t-r)|^2 \quad \text{(A.13)}$$

$$\leq -\frac{\theta q}{2|Q|}V(e(t)) + \frac{4}{\theta q}|Qp|^2|\eta(t)|^2 + \frac{4G^2}{\theta^3 q}|Q|^2|d(t-r)|^2$$

where $\eta(t) = x_1(t-r) - w(t)$, $\tilde{p}(x,e) = f(x+e) - f(x)$. Let $\sigma > 0$ be sufficiently small so that $4|Qp|(L+\theta)\exp(\sigma T_1)T_1\sqrt{\frac{|Q|}{a}} < q$ and $\sigma \leq \frac{\theta q}{8|Q|}$. The existence of sufficiently small $\sigma > 0$ satisfying the inequality $4|Qp|(L+\theta)\exp(\sigma T_1)T_1\sqrt{\frac{|Q|}{a}} < q$ is guaranteed by (3.17). Using (A.13), we conclude that:

$$V(t) \leq \exp(-4\sigma(t-r))V(r) + \frac{16|Q|}{\theta^2 q^2}|Qp|^2 \sup_{r\leq s\leq t}\left(\exp(-2\sigma(t-s))|\eta(s)|^2\right)$$

$$+ \frac{16G^2}{\theta^4 q^2}|Q|^3 \exp(2\sigma r) \sup_{0\leq s\leq t-r}\left(\exp(-2\sigma(t-s))|d(s)|^2\right) \quad \text{(A.14)}$$

for all $t \geq r$. Therefore, the following inequalities hold for all $t \geq r$:

$$|z(t)-x(t-r)| \le \exp(-2\sigma(t-r))\theta^{n-1}\sqrt{\frac{|Q|}{a}}|z(r)-x(0)| + \frac{4|Qp|}{q}\theta^{n-1}\sqrt{\frac{|Q|}{a}}\sup_{r\le s\le t}(\exp(-\sigma(t-s))|\eta(s)|)$$
$$+\frac{4|Q|}{q}\theta^{n-2}G\sqrt{\frac{|Q|}{a}}\exp(\sigma r)\sup_{0\le s\le t-r}(\exp(-\sigma(t-s))|d(s)|) \quad (A.15)$$

$$|z_i(t)-x_i(t-r)| \le \exp(-2\sigma(t-r))\theta^{i-1}\sqrt{\frac{|Q|}{a}}|z(r)-x(0)| + \frac{4|Qp|}{q}\theta^{i-1}\sqrt{\frac{|Q|}{a}}\sup_{r\le s\le t}(\exp(-\sigma(t-s))|\eta(s)|)$$
$$+\frac{4|Q|}{q}\theta^{i-2}G\sqrt{\frac{|Q|}{a}}\exp(\sigma r)\sup_{0\le s\le t-r}(\exp(-\sigma(t-s))|d(s)|) \quad (A.16)$$

where $a>0$ is a constant satisfying $a|x|^2 \le x'Qx$ for all $x\in\Re^n$. Using (3.2) and (A.16), we obtain for almost all $t \ge r$:

$$|\dot{w}(t)-\dot{x}_1(t-r)| \le \exp(-2\sigma(t-r))(L+\theta)\sqrt{\frac{|Q|}{a}}|z(r)-x(0)| + \frac{4|Qp|}{q}(L+\theta)\sqrt{\frac{|Q|}{a}}\sup_{r\le s\le t}(\exp(-\sigma(t-s))|\eta(s)|)$$
$$+\frac{4|Q|}{\theta q}G(L+\theta)\sqrt{\frac{|Q|}{a}}\exp(\sigma r)\sup_{0\le s\le t-r}(\exp(-\sigma(t-s))|d(s)|) + G|d_1(t-r)|$$

The above inequality implies that the following estimate holds for all $t\in[\tau_i,\tau_{i+1})$, where $\tau_i$ with $i\ge 1$ is an arbitrary sampling time with $\tau_i \ge r$:

$$|\eta(t)| \le |\xi(\tau_i)| + \exp(-\sigma(\tau_i-r))T_1(L+\theta)\sqrt{\frac{|Q|}{a}}|z(r)-x(0)|$$
$$+\frac{4|Qp|}{q}(L+\theta)T_1\sqrt{\frac{|Q|}{a}}\sup_{r\le s\le t}(\exp(-\sigma(\tau_i-s))|\eta(s)|)$$
$$+\frac{4|Q|}{\theta q}G(L+\theta)T_1\sqrt{\frac{|Q|}{a}}\exp(\sigma r)\sup_{0\le s\le t-r}(\exp(-\sigma(\tau_i-s))|d(s)|) + T_1 G\sup_{\tau_i-r\le s\le t-r}(|d_1(s)|)$$

Using the fact that $\tau_i \ge t-T_1$, in conjunction with the following inequalities

$$|\xi(\tau_i)| \le \exp(\sigma t)\sup_{\tau_i\le s\le t}(\exp(-\sigma(t-s))|\xi(s)|\exp(-\sigma s))$$
$$\le \exp(\sigma(t-\tau_i))\sup_{\tau_i\le s\le t}(\exp(-\sigma(t-s))|\xi(s)|)$$
$$\le \exp(\sigma T_1)\sup_{\tau_i\le s\le t}(\exp(-\sigma(t-s))|\xi(s)|)$$
$$\le \exp(\sigma T_1)\sup_{0\le s\le t}(\exp(-\sigma(t-s))|\xi(s)|)$$

$$\sup_{\tau_i-r\le s\le t-r}(|d_1(s)|) \le \exp(\sigma t)\sup_{\tau_i-r\le s\le t-r}(\exp(-\sigma(t-s))|d_1(s)|\exp(-\sigma s))$$
$$\le \exp(\sigma(t-\tau_i+r))\sup_{\tau_i-r\le s\le t-r}(\exp(-\sigma(t-s))|d_1(s)|)$$
$$\le \exp(\sigma(T_1+r))\sup_{\tau_i-r\le s\le t-r}(\exp(-\sigma(t-s))|d_1(s)|)$$
$$\le \exp(\sigma(T_1+r))\sup_{0\le s\le t}(\exp(-\sigma(t-s))|d_1(s)|)$$
$$\le \exp(\sigma(T_1+r))\sup_{0\le s\le t}(\exp(-\sigma(t-s))|d(s)|)$$

the above inequalities give for all $t\in[\tau_i,\tau_{i+1})$, where $\tau_i$ with $i\ge 1$ is an arbitrary sampling time with $\tau_i \ge r$:

$$|\eta(t)| \leq \exp(\sigma T_1) \sup_{0 \leq s \leq t} \left(\exp(-\sigma(t-s))|\xi(s)|\right) + \exp(-\sigma(t-r))T_1(L+\theta)\exp(\sigma T_1)\sqrt{\frac{|Q|}{a}}|z(r)-x(0)|$$
$$+ \frac{4|Q|p}{q}(L+\theta)\exp(\sigma T_1)T_1\sqrt{\frac{|Q|}{a}} \sup_{r \leq s \leq t}\left(\exp(-\sigma(t-s))|\eta(s)|\right) \quad \text{(A.17)}$$
$$+ T_1 G \exp(\sigma(r+T_1))\left(\frac{4|Q|}{\theta q}(L+\theta)\sqrt{\frac{|Q|}{a}} + 1\right) \sup_{0 \leq s \leq t}\left(\exp(-\sigma(t-s))|d(s)|\right)$$

Notice that the above inequality holds for all $t \geq r + T_1$. Setting $M := (L+\theta)\sqrt{\frac{|Q|}{a}}$, it follows from (A.17) and the inequality $4|Q|p(L+\theta)\exp(\sigma T_1)T_1\sqrt{\frac{|Q|}{a}} < q$ that the following inequality holds for all $t \geq r + T_1$:

$$\sup_{r+T_1 \leq s \leq t}\left(\exp(\sigma s)|\eta(s)|\right) \leq \frac{q \exp(\sigma T_1)}{q - 4|Q|pM\exp(\sigma T_1)T_1} \sup_{0 \leq s \leq t}\left(\exp(\sigma s)|\xi(s)|\right) + \frac{qM \exp(\sigma(r+T_1))T_1}{q-4|Q|pM\exp(\sigma T_1)T_1}|z(r)-x(0)|$$
$$+ \sup_{r \leq s \leq r+T_1}\left(\exp(\sigma s)|\eta(s)|\right) + \frac{(4|Q|M + q\theta)\exp(\sigma(r+T_1))}{\theta(q-4|Q|pM\exp(\sigma T_1)T_1)} GT_1 \sup_{0 \leq s \leq t}\left(\exp(\sigma s)|d(s)|\right) \quad \text{(A.18)}$$

The existence of constants $\sigma > 0$, $A_i > 0$ ($i=1,\ldots,4$), which are independent of $T_2 > 0$ and $l,m$, satisfying (3.23) is a direct consequence of (A.15) and the above inequality. The proof is complete. ◁

**Proof of Lemma 3.6:** Let $\sigma > 0$ be sufficiently small such that (3.23) holds and such that $|k|T_2 + C|k|\exp(\sigma(T_2 + r + \tau)) < 1$, where $C := K\frac{((nL+1)T)^{l+1}}{1-(nL+1)T}$. The existence of sufficiently small $\sigma > 0$ satisfying $|k|T_2 + C|k|\exp(\sigma(T_2 + r + \tau)) < 1$ is guaranteed by (3.19). Using (3.16), we obtain for all $i \in Z^+$ and $t \in [iT_2 + \tau, (i+1)T_2 + \tau)$:

$$|u(t-\tau) - k'x(t)| \leq |k|\left|\Phi_{l,m}(z(iT_2), \breve{T}_{r+\tau}(iT_2)u) - x(t)\right|$$
$$\leq |k|\left|\Phi_{l,m}(z(iT_2), \breve{T}_{r+\tau}(iT_2)u) - x(iT_2 + \tau)\right| + |k||x(iT_2 + \tau) - x(t)| \quad \text{(A.19)}$$

Using (3.13), we obtain for all $i \in Z^+$ with $iT_2 \geq r$ and $t \in [iT_2 + \tau, (i+1)T_2 + \tau)$:

$$\left|\Phi_{l,m}(z(iT_2), \breve{T}_{r+\tau}(iT_2)u) - x(iT_2 + \tau)\right| \leq K\frac{((nL+1)T)^{l+1}}{1-(nL+1)T}\left(|z(iT_2)| + \sup_{iT_2-r-\tau \leq s < iT_2}|u(s)|\right)$$
$$+ \exp((nL+1)(r+\tau))\left(G(r+\tau) \sup_{iT_2-r \leq s \leq iT_2+\tau}|d(s)| + |z(iT_2) - x(iT_2-r)|\right) \quad \text{(A.20)}$$

Combining (A.19) and (A.20) we obtain for all $i \in Z^+$ with $iT_2 \geq r$ and $t \in [iT_2 + \tau, (i+1)T_2 + \tau)$:

$$|u(t-\tau) - k'x(t)| \leq K|k|\frac{((nL+1)T)^{l+1}}{1-(nL+1)T} \sup_{iT_2-r \leq s < iT_2+\tau}|u(s-\tau) - k'x(s)|$$
$$+ |k|\left(K\frac{((nL+1)T)^{l+1}}{1-(nL+1)T} + \exp((nL+1)(r+\tau))\right)|z(iT_2) - x(iT_2-r)|$$
$$+ G(r+\tau)|k|\exp((nL+1)(r+\tau)) \sup_{iT_2-r \leq s \leq iT_2+\tau}|d(s)| \quad \text{(A.21)}$$
$$+ K|k|(1+|k|)\frac{((nL+1)T)^{l+1}}{1-(nL+1)T} \sup_{iT_2-r \leq s \leq iT_2+\tau}|x(s)| + |k||x(iT_2+\tau) - x(t)|$$

On the other hand, using (3.2) and (3.3), we conclude that the following inequality holds for all $i \in Z^+$ and $t \in [iT_2 + \tau, (i+1)T_2 + \tau)$:

$$\exp(\sigma t)|x(t) - x(iT_2 + \tau)| \leq T_2(nL+1)\exp(\sigma t)\sup_{iT_2+\tau \leq s \leq t}|x(s)| + T_2 \exp(\sigma t)|u(t-\tau)| + T_2 G \exp(\sigma t)\sup_{iT_2+\tau \leq s \leq t}|d(s)|$$

$$\leq T_2(nL+1)\exp(\sigma t)\sup_{iT_2+\tau \leq s \leq t}|x(s)| + T_2 \exp(\sigma t)|u(t-\tau) - k'x(t)| + T_2|k|\exp(\sigma t)|x(t)| + T_2 G \exp(\sigma t)\sup_{iT_2+\tau \leq s \leq t}|d(s)|$$

$$\leq T_2(nL+1+|k|)\exp(\sigma T_2)\sup_{iT_2+\tau \leq s \leq t}(\exp(\sigma s)|x(s)|) + T_2 \exp(\sigma t)|u(t-\tau) - k'x(t)| + T_2 G \exp(\sigma T_2)\sup_{0 \leq s \leq t}(\exp(\sigma s)|d(s)|)$$

(A.22)

Inequality (A.21) implies that the following inequality holds for all $i \in Z^+$ with $iT_2 \geq r$ and $t \in [iT_2 + \tau, (i+1)T_2 + \tau)$:

$$\exp(\sigma t)|u(t-\tau) - k'x(t)| \leq K|k|\frac{((nL+1)T)^{l+1}}{1-(nL+1)T}\exp(\sigma(T_2 + r + \tau))\sup_{iT_2 - r \leq s < iT_2 + \tau}(\exp(\sigma s)|u(s-\tau) - k'x(s)|)$$

$$+ \exp(\sigma t)|k|\left(K\frac{((nL+1)T)^{l+1}}{1-(nL+1)T} + \exp((nL+1)(r+\tau))\right)|z(iT_2) - x(iT_2 - r)|$$

$$+ G(r+\tau)|k|\exp((nL+1)(r+\tau))\exp(\sigma(T_2 + r + \tau))\sup_{iT_2 - r \leq s < iT_2 + \tau}(\exp(\sigma s)|d(s)|)$$

$$+ K|k|(1+|k|)\frac{((nL+1)T)^{l+1}}{1-(nL+1)T}\exp(\sigma(T_2 + r + \tau))\sup_{iT_2 - r \leq s < iT_2 + \tau}(\exp(\sigma s)|x(s)|) + \exp(\sigma t)|k||x(iT_2 + \tau) - x(t)|$$

(A.23)

It follows from Lemma 3.5 and inequality (3.23) that the following inequality holds for all $i \in Z^+$ with $iT_2 \geq r + T_1$ and $t \in [iT_2 + \tau, (i+1)T_2 + \tau)$:

$$\exp(\sigma t)|z(iT_2) - x(iT_2 - r)| \leq A_1 \exp(\sigma(T_2 + \tau + r))|z(r) - x(0)| + A_2 \exp(\sigma(T_2 + \tau))\sup_{0 \leq s \leq t}(\exp(\sigma s)|\xi(s)|)$$

$$+ A_3 \exp(\sigma(T_2 + \tau))\sup_{r \leq s \leq r + T_1}(\exp(\sigma s)|w(s) - x_1(s-r)|) + A_4 \exp(\sigma(T_2 + \tau))\sup_{0 \leq s \leq t}(\exp(\sigma s)|d(s)|)$$

(A.24)

Combining (A.23) and (A.24) we obtain for all $i \in Z^+$ with $iT_2 \geq r + T_1$ and $t \in [iT_2 + \tau, (i+1)T_2 + \tau)$:

$$\exp(\sigma t)(1-|k|T_2)|u(t-\tau) - k'x(t)| \leq C|k|\exp(\sigma(T_2 + r + \tau))\sup_{iT_2 - r \leq s < iT_2 + \tau}(\exp(\sigma s)|u(s-\tau) - k'x(s)|)$$

$$+ A_1 \exp(\sigma(T_2 + r + \tau))|k|(C + \exp((nL+1)(r+\tau)))|z(r) - x(0)|$$

$$+ A_2 \exp(\sigma(T_2 + \tau))|k|(C + \exp((nL+1)(r+\tau)))\sup_{0 \leq s \leq t}(\exp(\sigma s)|\xi(s)|)$$

$$+ A_3 \exp(\sigma(T_2 + \tau))|k|(C + \exp((nL+1)(r+\tau)))\sup_{r \leq s \leq r + T_1}(\exp(\sigma s)|w(s) - x_1(s-r)|)$$

$$+ |k|\exp(\sigma(T_2 + \tau))[A_4 C + \exp((nL+1)(r+\tau))(A_4 + G(r+\tau)\exp(\sigma r)) + T_2 G \exp(-\sigma\tau)]\sup_{0 \leq s \leq t}(\exp(\sigma s)|d(s)|)$$

$$+ |k|\exp(\sigma T_2)[C(1+|k|)\exp(\sigma(r+\tau)) + T_2(nL+1+|k|)]\sup_{iT_2 - r \leq s \leq t}(\exp(\sigma s)|x(s)|)$$

Inequality (3.24) is a direct consequence of the above inequality. Indeed, using the above inequality and inequality $|k|T_2 + C|k|\exp(\sigma(T_2 + r + \tau)) < 1$ we can compute an upper bound for $\sup_{jT_2 + \tau \leq s \leq t}(\exp(\sigma s)|u(s-\tau) - k'x(s)|)$, where $j = \min\{j \in Z^+ : jT_2 \geq r + T_1\}$. The proof is complete. ◁

31